\documentclass[12pt,twoside]{amsart}

\usepackage{amsfonts}
\usepackage{amsmath}
\usepackage{latexsym}

\newtheorem{thm}{{\sc Theorem}}[section]
\newtheorem{prop}[thm]{{\sc Proposition}}
\newtheorem{lem}[thm]{{\sc Lemma}}
\newtheorem{cor}[thm]{{\sc Corollary}}
\newtheorem{defn}[thm]{{\sc Definition}}
\newtheorem{claim}[thm]{{\sc Claim}}
\newtheorem{question}[thm]{{\sc Question}}

\newcommand{\qqed}{\hspace*{\fill} $\Box$}

\title[Higher direct images]{Higher direct images of dualizing sheaves
                            of Lagrangian fibrations}
\author{Daisuke Matsushita}
\address{Research Institute for Mathematical Sciences \\
         Kyoto University, Oiwake-Cyo Kitashirakawa \\
         Sakyo-Ku Kyoto 606-8052 Japan}
\thanks{*Research Fellow of the Japan Society for the Promotion of Science} 
\email{tyler@kurims.kyoto-u.ac.jp}

\begin{document}
\maketitle
\begin{abstract}
 We prove that 
 the higher direct images of the dualizing sheaf of 
 a Lagrangian fibration 
 between smooth projective manifolds 
 are isomorphic to the cotangent bundles of base space.
 As a corollary, we obtain that every Hodge number 
 of the base space of a fibre space of an irreducible symplectic
 manifold is the same to that of a projective space
 if the base space is smooth.
\end{abstract}
\section{Introduction} 
 We begin with the definition of {\it  Lagrangian fibrations\/}.
\begin{defn}
 Let $X$ be a K\"{a}hler manifold with a
 holomorphic symplectic form $\omega$
 and $S$ a normal variety. A proper surjective morphism $f : X \to S$
 is said to be a Lagrangian fibration if 
 a general fibre $F$ of $f$ is a Lagrangian submanifold with
 respect to $\omega$,
 that is,
 the restriction of 2-form $\omega |_{F}$ is 
 identically zero and $\dim F = (1/2)\dim X$.
\end{defn}

\noindent
 In this note, we investigate higher direct images of the
 dualizing sheaf of a Lagrangian fibration. Our result
 is the following.
\begin{thm}\label{main}
 Let $f :X \to S$ be a Lagrangian fibration between 
 smooth projective manifolds. Then
$$
 R^{i}f_{*}{\mathcal O}_{X} \cong \Omega^{i}_{S}.
$$
 In particular, $R^{i}f_{*}{\mathcal O}_{X}$ are locally free.
\end{thm}

\noindent
{\sc Remark 1. \quad} It is known that the higher direct images 
 of the dualizing sheaf are locally free if the discriminant
 locus is a normal crossing divisor
 \cite[Theorem 2.6]{torsionfree}, 
 \cite[Theorem 2.6]{moriwaki},
 and \cite[Theorem 1]{nakayama}. 
 However, there
 is a projective Lagrangian fibration
 whose discriminant locus is not a normal crossing divisor.
 The above theorem
 asserts that the higher direct images are locally free even if
 the discriminant locus is not necessarily normal crossing.
 
\vspace{3mm}

\noindent
{\sc Example. \quad} Let $X := {\rm Hilb}^{2}T$, where $T$ is
 a projective $K3$ surface 
 which has an elliptic fibration $T \to {\mathbb P}^{1}$.
 Then $X$ admits a Lagrangian fibration over ${\mathbb P}^{2}$ and
 whose discriminant locus is not a normal crossing
 divisor.

\vspace{3mm}

\noindent
{\sc Remark 2. \quad} The assumption  projectivity 
 is crucial in our present proof. The author suspect,
 however, that
 the most part of our results holds true for  Lagrangian
 fibrations. Please see {\sc Remark 4} in Section 3.

\vspace{3mm}

\noindent 
 Combining Theorem \ref{main} with
 \cite[Theorem 2]{matsu} and \cite[Theorem 1]{matsu2},
 we obtain the following corollary.
\begin{cor}\label{main_result}
 Let $f :X \to S$ be a fibre space of a projective irreducible
 symplectic manifold over a smooth projective base $S$. 
 Then $S$ is Fano and
 hodge numbers of $S$
 are the following:
$$
 h^{q}(S, \Omega^{p}_{S})= 
\left\{
\begin{array}{c}
  1 \quad ( p = q )\\
  0 \quad ( p \ne q ).
\end{array}
\right.
$$
 Every Hodge number of $S$ is the same as that of
 ${\mathbb P}^{n}$, where $n = \dim S$.
\end{cor}

\noindent
{\sc Remark 3. \quad} 
 Every Hodge number of a
  quadric hypersurface $Q$ in ${\mathbb P}^{2n}$
 is  same to that of ${\mathbb P}^{2n-1}$ and $Q$ is
 not isomorphic to ${\mathbb P}^{2n-1}$ $( n \ge 2)$. The author
 does not know whether there exists an example such that
 a fibre space of an irreducible symplectic manifold $f: X \to S$
 whose base space $S$ is not isomorphic to a projective space.
 Recently, Miyaoka \cite{miyaoka} announced that $S$ is isomorphic to
 a projective space if $f$ admits a local section. 
 Note that if $f$ admits a
 local section, then $S$ is smooth. The author would like to
 state one question.
\begin{question}\label{ampleness}
 Under the assumptions of
 Theorem \ref{main}, we have 
$$
 R^{n-1}f_{*}\omega_{X/S} \cong T_{S},
$$
 which is the dual of $R^{1}f_{*}{\mathcal O} \cong \Omega^{1}_{S}$.
 Note that
 $\omega_{X/S}$ is the relative dualizing sheaf
 and $n = \dim S$.
 Then $R^{n-1}f_{*}\omega_{X/S}$ is nef? In particular,
 $R^{n-1}f_{*}{\omega_{X/S}}$ is ample if $X$ is
 an irreducible symplectic manifold?
\end{question}
 \noindent
 Note that it is known that $R^{i}f_{*}\omega_{X/S}$ is weakly
 positive.
\begin{cor}\label{4dim}
 If $\dim X =4$, then $T_S$ is nef. In particular, $T_S$ is
 ample if $X$ is an irreducible symplectic manifold. 
\end{cor}
\noindent
  This paper is organized as follows. In section 2, we prepare for
 the proof of Theorem \ref{main}. Theorem \ref{main},
 Corollary \ref{main_result} 
 and Corollary \ref{4dim} are
 proved in section 3. Section 4 is devoted to the 
 classification of singular fibres over codimension one
 points which needs for the proof 
 of Theorem \ref{main}.

\vspace{3mm}
\noindent
{\sc Acknowledgement. \quad}  The author would like to express
 his gratitude to Professors A.~Fujiki, Y.~Miyaoka, S.~Mori and
 N.~Nakayama for
 their encouragement and suggestions. 
 The author also thanks to Dr.~O.~Fujino who gave
 him useful comments. The main idea of this note is produced
 in author's stay at \'{E}cole Normale Sup\'{e}rieure, Paris.
 The author wishes to thank to Professor
 A.~Beauville who gave him an opportunity to stay at E. N. S.

\section{Preliminary}
 We correct definitions and fundamental material which are
 needed in the proof of Theorem \ref{main}.
 First we consider a relative compactification
 of an Abelian fibration.
\begin{prop}\label{toroidal_model}
 Let $f^{\circ} : X^{\circ} \to \Delta^{*}(t_1)
 \times \Delta^{n-1}(t_2 , \cdots , t_n)$ 
 be a smooth projective Abelian
 fibration. Assume that
 $f^{\circ}$ has the following properties:
\begin{enumerate}
 \item $f^{\circ}$ has a projective relative compactification
       $f' : Y \to \Delta^{n}$.
 \item $f^{' -1}(p)$ has a reduced component for  
       a general point $p$ of the discriminat locus of $f'$.
 \item The monodoromy matrix $T$ of $R^{1}f^{\circ}_{*}{\mathbb C}$ 
       around $\Delta^{*}(t_1)$ is unipotent.
 \item The rank of $T-I$ is at most one, where $I$ is the identity
       matrix.
\end{enumerate}
 Then there exists a projective relative compactification
 $f : X \to \Delta^{n}$
 which has the following properties:
\begin{enumerate}
 \item Every singular fibre of $f$ is an Abelian variety
       or a cycle of several copies of 
       a ${\mathbb P}^{1}$-bundle over an Abelian
       variety. In later case, $E_i \cap E_j$ forms
       a section of the ruling of $E_i$ and $E_j$,
       where $E_i$ and $E_j$ are irreducible components
       of a singular fibre.
 \item The morphism $f^{-1}(D) \to D$ is a locally trivial
       deformation, where $D$ is the discriminat locus of $f$.
 \item $K_X$ is $f$-trivial.
 \item Every fibre of $f$ is reduced.
\end{enumerate}
\end{prop}
\noindent
{\sc Proof. \quad} This proposition is essentially deduced
 from \cite[Theorem 5.3]{namikawa}. However, for a complete sake
 and the proof of projectivity of the relative compactification 
 is not published, we give a proof. 
 Let $\tau (t_1, {\bf t})$ be the period matrix of $f^{\circ}$
 where ${\bf t} = (t_2 ,\cdots , t_n)$, 
 $m$ is the dimension of a fibre of $f^{\circ}$,
 $I_k$ is the identity matrix of rank $k$
 and $\delta$ is the polarization matrix of $f$.
 We consider $T$ as an element of ${\rm Sp}(m, {\mathbb R})$
 by the conjugation
$$
 T' := 
 \left(
\begin{array}{cc}
 I_m & 0  \\
 0 & \delta^{-1}
\end{array}
\right)
 T \left(
\begin{array}{cc}
 I_m &  0 \\
 0 & \delta
\end{array} 
\right).
$$
 Let $G$ be the semiproduct of ${\mathbb Z}$ and
 $\chi := \{a + b\delta ; a,b \in {\mathbb Z}^{\oplus m} \}$
 with commuting relations 
 $(a + b\delta )\cdot \gamma = \gamma \cdot (a+ b\delta)T' $
 for $\gamma \in {\mathbb Z}$ and 
 $a+ b\delta \in \chi$.
 By the assumption, $f'$ admits a meromorphic section locally over 
 an open set $V$ such that ${\rm codim} (\Delta^{n}\setminus V) \ge 2$.
 Since $\pi_{1}(V) \cong \pi_{1}(\Delta^{n}) = \{1\}$,
 $f^{\circ}$ admits a global section. Hence
 the original $X^{\circ}$ is reconstructed
 as follows:
$$
 X^{\circ} = 
 ({\mathbb H}\times \Delta^{n-1} \times {\mathbb C}^{m})/G
$$
 where ${\mathbb H}$ is an upper half space.
 The action of $G$ is defined by
$$
\begin{array}{cccc}
 \gamma : & {\mathbb H}\times \Delta^{n-1}\times {\mathbb C}^{m} &
           \to &  {\mathbb H}\times \Delta^{n-1}\times {\mathbb C}^{m} \\
     & (s , {\bf t}, {\bf z}) 
     & \mapsto & (s+1, t , {\bf z} 
     (C\tau(\exp (2\pi is),{\bf t}) + D)^{-1}) \\
 a + b\delta : &  {\mathbb H}\times \Delta^{n-1}\times {\mathbb C}^{m}
        & \to &  {\mathbb H}\times \Delta^{n-1}\times {\mathbb C}^{m} \\
      & (s , {\bf t}, {\bf z}) 
      & \mapsto & (s , t, {\bf z} + (a + b\delta) 
\left(
 \begin{array}{c}
  \tau (\exp (2\pi i s),{\bf t}) \\
   I_m
 \end{array}
\right)), 
\end{array}
$$
 where
$$
 T' = 
\left(
\begin{array}{cc}
 A & B \\
 C & D
\end{array}
\right).
$$
\begin{lem}\label{toroidal_model_I}
 If the rank of $T - I_{2m}$ is zero, 
 then there exists a  smooth Abelian fibration
 $f : X \to \Delta^{n}$ which is a relative compactification
 of $X^{\circ}$. This fibration has properties of
 Proposition \ref{toroidal_model} (1), (2), (3) and (4). 
\end{lem} 
\noindent
{\sc Proof. \quad} By the assumption,
 the monodoromy matrix $T = I_{2m}$. 
 Therefore the period $\tau$ can be extended
 holomorphically over $\Delta^{n}$. We construct $X$ as follows:
$$
 X = (\Delta^{n}\times {\mathbb C}^{m})/\chi,
$$
 where 
$ 
 \chi := \{ a + b\delta ; a,b \in {\mathbb Z}^{\oplus 2m}\}
$
 and the action of $\chi$ is
$$
\begin{array}{cccc}
 a + b\delta : &  \Delta^{n}\times {\mathbb C}^{m}
        & \to &  \Delta^{n}\times {\mathbb C}^{m} \\
      & (t_1 ,{\bf t}, {\bf z}) & \mapsto
 & (t_1 ,{\bf t}, {\bf z} + (a + b\delta) 
\left(
 \begin{array}{c}
  \tau (t_1 ,{\bf t}) \\
   I_m
 \end{array}
\right)).
\end{array} 
$$
 Since $f : X \to \Delta^{n}$ is a smooth morphism,
 it is obvious that $X$ satisfies the assertions of Lemma.
\qed

\noindent
 We consider the case that the rank of $T - I_{2m}$ is one.
 Since the discriminant locus of $g$ is smooth,
 $\lim_{t_1 \to 0} \tau(t_1 , {\bf t})$ exists in a rational boundary 
 component $F_{\alpha}$ of Siegel upper half plane
 and $F_{\alpha}$ does not depend ${\bf t}$.
 (\cite[Remark (4.4)]{namikawa}). 
 There exists a matrix $M \in {\rm Sp}(m,{\mathbb Z})$ 
 which transform $F_{\alpha}$ to $F_{g'}$
 (cf. \cite[(2.3)]{namikawa}),
 where $F_{g'}$ is the rational boundary component of 
 Siegel upper half plane defined by
$$
 F_{g'} := \left\{
\left(
\begin{array}{cc}
 0 & 0 \\
 0 & Z''
\end{array}
\right) ;
 Z'' \in {\rm GL}(m-g' ,{\mathbb C}) , I_{m-g'} - Z''\bar{Z''} > 0
\right\}.
$$
 Since $T'$ is unipotent,
 the transformation of $T'$ by $M$ can be written
\begin{equation}\label{Monodoromy_Matrix_canonical_form}
 MT'M^{-1} = 
\left( 
\begin{array}{cc}
 I_m & B \\
 0   & I_m
\end{array}
\right), B = \left(
 \begin{array}{cc}
  0 & 0 \\
  0 & b
 \end{array}
\right),
\end{equation}
 by \cite[(2.6)]{namikawa}.
 Note that $b \in {\mathbb R}^{1}$ because
 the rank of $T - I_{2m}$ is one.
 From this representation 
 of the monodoromy matrix $MT' M^{-1}$,
 the period matrix can be written as
$$
 M \tau (t_1 ,{\bf t} ) = 
\left(
\begin{array}{cc}
 \tau' (t_1 ,{\bf t}) & \tau''' (t_1 ,{\bf t}) \\
 {}^{t}\tau''' (t_1 ,{\bf t}) & \tau'' (t_1 ,{\bf t})
\end{array}
\right)
 +
 \frac{1}{2\pi i} 
\left(
\begin{array}{cc}
0  & 0 \\
0  & b
\end{array}
\right) \log t_1 ,
$$
 where each $\tau^{*}(t_1 ,{\bf t})$ be a holomorphic function on
 $\Delta^{n}$. From \cite[(2.9)]{namikawa},
 ${\rm Im}\tau' > 0$ and ${\rm Im}\tau'' > 0$.
 We denote each element of $\chi M^{-1}$ as 
$$
 \{(\alpha_1 , \alpha_2 , \beta_1  , \beta_2 ) \in \chi M^{-1};
   \alpha_1 , \beta_1  \in {\mathbb R}^{m-1}, \alpha_2 , \beta_2 
  \in{\mathbb R}^{1}\}.
$$
 Let $\bar{\chi}$ be the
 sublattice of 
 $\chi M^{-1}$ such that
$$
 \bar{\chi} := \{ (\alpha_1 , \alpha_2 , \beta_1  , \beta_2) \in \chi M^{-1}
             ; \alpha_1 = \alpha_2 = \beta_1  = 0 \}.
$$
 and $(0,0,0,c)$ its generator.
 We define a toric variety ${\mathcal M}$ for the construction $X$.
\begin{defn}
 Let $\sigma_{k}$ be a cone defined by
$$ 
 \sigma_{k} := \{ (u, y) \in 
                  {\mathbb R}^{2} ; k u \le y \le (k+1)u \} 
$$
 Then ${\mathbb C}[\sigma_{k}^{\vee}\cap {\mathbb Z}^{2}]$ 
 is isomorphic to ${\mathbb C}$-subalgebra
 of ${\mathbb C}[v^{\pm 1}, t_{1}^{\pm 1}]$ generated by
 monomials $\{ t_{1}^{-k}v, t_{1}^{k+1}v^{-1} \}$ and
 there is a natural morphism
$
 {\rm Spec}{\mathbb C}[\sigma_{k}^{\vee} \cap {\mathbb Z}^{2} ]
 \to
 {\rm Spec}{\mathbb C}[t_1]$. We define
 the toric variety ${\mathcal M}$ 
 as follows:
$$
 {\mathcal M} := \bigcup_{k\in {\mathbb Z}} \left(
                         {\rm Spec}{\mathbb C}[\sigma_{k}^{\vee}
                          \cap {\mathbb Z}^{2}] 
                         \right)
                         \times_{{\mathbb C}^{1}}
                         \Delta^{1}.
$$
 Note that ${\mathcal M}$ has a fibration ${\mathcal M} \to \Delta^{1}$
 and whose singular fibre is an infinite chain of ${\mathbb P}^{1}$.
\end{defn}
\begin{lem}\label{logarithmic_version}
 Let ${\mathcal M}^{*} :=
 {\mathcal M}
 \times_{\Delta^{1}}\Delta^{*}$.
 Then ${\mathcal M}^{*} \cong \Delta^{*}\times {\mathbb C}^{*}$
 and the variable $v$ is considered to be a coordinate of
 ${\mathbb C}^{*}$. 
 We define the action of $\chi M^{-1}$
 on $\Delta^{n-1} \times
 {\mathbb C}^{m-1} 
 \times {\mathcal M}^{*}$
 as follows:
$$
 \Gamma_{(\alpha_1 , \alpha_2 , \beta_1  ,\beta_2 )} 
     : ({\bf t}, t_1 ,{\bf z},v) \mapsto 
        ({\bf t}, t_1 , {\bf z} + A , Bv ),
$$
 where 
\begin{eqnarray*}
       A &=& \alpha_1 \tau'(t_1 ,{\bf t}) 
       + \beta_1  + \alpha_2 \tau'''(t_1 ,{\bf t}) \\ 
       B &=& \exp \left(\frac{2\pi i}{c} 
        (\alpha_1 {}^{t}\tau'''(t_1 ,{\bf t}) 
        + \alpha_2 \tau''(t_1 ,{\bf t}) + \beta_2) \right)
         {t_1}^{\alpha_2 b /c }.
\end{eqnarray*}
  Then the quotient 
 $(\Delta^{n-1}\times {\mathcal M}^{*} \times 
  {\mathbb C}^{m-1})/\chi M^{-1} \cong X^{\circ}$.
\end{lem}
\noindent
{\sc Proof. \quad}
 We bigin with to prove that the action is well defined.
 It is enough
 to show that $\alpha_2 b /c$ is an integer
 for every element of $\chi M^{-1}$.
 For each element $x$ of $\chi$, 
$
 xT' \in \chi .
$
 Hence
$ 
 \chi (T' - I_{2m}) \subset \chi 
$
 and we have
$
 \chi M^{-1}(MT' M^{-1} -  I_{2m} ) \subset \chi M^{-1}.
$
 On the contrary,
$
 (0,0,0, \alpha_2 b) = (\alpha_1 , \alpha_2 ,\beta_1  , \beta_2 )
 (M T'M^{-1} - I_{2m}).
$
 Therefore
 $(0,0,0, \alpha_2b) \in \bar{\chi}$ for
 every element of $\chi M^{-1}$ and
 $\alpha_2 b/c$ is an integer. 
 If we  consider the morphism 
\begin{eqnarray*}
  \Delta^{n-1}\times \Delta^{*}\times
  {\mathbb C}^{m-1} \times  {\mathbb C}
  & \to &
  \Delta^{n-1}\times \Delta^{*}\times 
  {\mathbb C}^{m-1}  \times  {\mathbb C}^{*}   \\
 ({\bf t}, t_1 ,  {\bf z},z) & \mapsto & 
 ({\bf t}, t_1 , {\bf z}, \exp(2\pi i z/c)),
\end{eqnarray*}
 and
 rewrite the construction of $X^{\circ}$ by
 $MT' M^{-1}$ and $M\tau$,
 we obtain that 
 $
 \Delta^{n-1}\times {\mathcal M}^{*} \times {\mathbb C}^{m-1}/
 \chi  M^{-1} \cong X^{\circ}.
 $
\qed

\begin{lem}\label{toroidal_model_II}
 The action of $\chi M^{-1}$ in 
 Lemma \ref{logarithmic_version} can be
 extended on $\Delta^{n-1}\times {\mathbb C}^{n-1} \times {\mathcal M}$.
 This action 
 is properly
 discontinuous and fixed points free. 
 The quotient $X$ is a relative compactification of $X^{\circ}$ and
 $f : X \to \Delta^{n}$ 
 has properties of Proposition \ref{toroidal_model} (1), (2), (3) and (4).
\end{lem} 
\noindent
{\sc Proof. \quad} 
 By definition, every singular fibre of 
 $\Delta^{n-1}\times {\mathcal M} \times {\mathbb C}^{n-1} 
  \to \Delta^{n}$ is the product of an infinite chain of
 ${\mathbb P}^{1}$ and ${\mathbb C}^{n-1}$.
 We denote each irreducible component of
 a singular fibre by $V_k$. 
 Since the action of $\chi M^{-1}$ on ${\mathcal M}$
 maps ${\rm Spec} {\mathbb C}[\sigma_{k}\cap {\mathbb Z}^2]$
 to  ${\rm Spec} {\mathbb C}
 [\sigma_{k + (\alpha_2 b)/c}\cap {\mathbb Z}^2]$, 
 $V_k$ maps to $V_{k + \alpha_{2} b/c}$ by this action. 
 Hence this action
 is  properly discontinuous
 and fixed point free. 
 The quotient $X$ is smooth and a relative compactification
 of $X^{\circ}$. By construction, 
 it is obvious that $f$ satisfies the properties
 (1) and (2) of Propositon \ref{toroidal_model}.
 We will prove that every fibre of $f$ is
 reduced. By construction, the multiplicity of
 each component of a singular fibre is same. Thus
 it is enough to prove that every fibre has a reduced
 component.
 The section 
 $\Delta^{n} \ni ({\bf t},t_1 ) \mapsto ({\bf t},t_1  , 1, 0)
 \in \Delta^{n-1}\times {\mathcal M}\times {\mathbb C}^{n-1}$
 induces a global section of $f$. Hence
 every fibre of $f$ is reduced.
 Moreover, $K_X$ is $f$-trivial because there is 
 the holomorphic $n+m$-form
$$ 
 dt_{1} \wedge \cdots \wedge dt_{n} \wedge 
 dz_1 \wedge \cdots \wedge dz_{m-1} \wedge \frac{dv}{v}
$$
 which is  $\chi M^{-1}$-invariant and nowhere vanishing
 on 
 $\Delta^{n-1}\times {\mathcal M}\times {\mathbb C}^{m-1}$.
\qed

\noindent
 The rest of assertions of Proposition \ref{toroidal_model} 
 is deduced from the following Lemma.
\begin{lem}\label{quotient}
 Let $f : X \to \Delta^{n}$ be the relative compactification
 constructed in Lemma \ref{toroidal_model_I} or Lemma
 \ref{toroidal_model_II}.
 Then $f$
 is projective.
\end{lem}
\noindent
{\sc Proof. \quad} 
 By changing birational model,
 we may assume that there is a morphism $Y \to X$. 
 Let $H$ be a relative ample divisor on $Y$ and $H'$ its proper transform
 on $X$. Then $H'$ is a big divisor on each component of a fibre of $f$.
 If $T = I_{2m}$, every fibre of $f$ is an Abelian variety.
 Thus every big divisor on each component
 of a fibre of $f$ is ample. Therefore $H'$ is ample.
 If $T \ne I_{2m}$, we need the following
       claim:
\begin{claim}\label{projectivity}
   Let $f : X \to \Delta^{n}$ be the relative compactification 
   constructed in Lemma \ref{toroidal_model_II} and
   $E_s$ an irreducible component of a singular fibre of $f$
   at $s$.
   Then there is a section $e$ of the ruling
   of $E_s$  such that $-K_{E_s} \equiv 2e$ and $e$ is nef.
\end{claim}
\noindent
{\sc Proof. \quad}
   Let $C$ be a smooth curve on $\Delta^{n}$ such that
   $C \cap D = \{s\}$, where $D$ is the
   discriminat locus of $f$. We consider the
   restriction 
   $f_{C} : X_{C}:= X\times_{\Delta^{n}}C \to C$. 
   We denote each irreducible component
   of the singular fibre of $f_{C}$ by $E_{i,s}$.  
   Let $E_{i+1 ,s}$ and $E_{i-1,s}$ are next
   components of $E_{i,s}$ and
   $F_{ij ,s} := E_{i,s} \cap E_{j,s}$.
   Since $K_{X_{C}}$ is $f_C$-trivial,
   $K_{E_{i,s}} \sim - F_{i i+1 ,s} - F_{i i-1 ,s}$.
   Each $E_{i,s}$ is a ${\mathbb P}^{1}$-bundle
   over an Abelian variety. Hence $-K_{E_{i,s}}$ is nef.
   In order to prove Claim, it is enough to show
   that $F_{i i+1 ,s}$ and $F_{i  i-1 ,s}$ are
   numerically equivalent on $E_{i ,s }$.
   Since $X_s$ is a deformation retract of $X_C$, there is
   a corrasping morphism $c_t : X_t \to X_s$, where
   $X_t$ is a general fibre of $f_{C}$.
   Let $E'_{i,s} := c^{-1}_{t}(E_{i,s})$ and
   $F'_{ij ,s} := c^{-1}_{t}(F_{ij, s})$.
   Note that $E'_{i,s} \to E_{i,s}$ is a real blow up
   along $F_{ij,s}$ and $F'_{ij,s}$ is a $S^{1}$-bundle
   over $F_{ij,s}$.
   We consider the following Mayer-Vietoris sequences:
\begin{eqnarray*}
 && 0 \to {\mathbb C}_{X_{t}} \to \oplus_{E'_{i ,s}} {\mathbb C}_{E'_{i,s}}
   \to \oplus_{F'_{ij,s}} {\mathbb C}_{F'_{ij,s}} \to 0 \\
 && 0 \to {\mathbb C}_{X_{s}} \to \oplus_{E_{i ,s}} {\mathbb C}_{E_{i,s}}
   \to \oplus_{F_{ij,s}} {\mathbb C}_{F_{ij,s}} \to 0
\end{eqnarray*}
   From the above sequence,
$$
 0 \to H^{0}(X_{t}) \to \oplus_{E'_{i,s}} H^{0}(E'_{i,s}) \to
 \oplus_{F'_{ij,s}} H^{0}(F'_{ij ,s}) 
 \stackrel{\alpha_1}{\to} H^{1}(X_{t}),
$$ 
   and
$$
\begin{array}{rcccccc}
   & & 0 & & 0 & & \\
  & & \downarrow & & \downarrow & &  \\        
 \to H^{1}(X_{s}) & \to & \oplus_{E_{i,s}} H^{1}(E_{i,s}) & 
     \stackrel{\alpha_4}{\to} &
     \oplus_{F_{ij,s}} H^{1}(F_{ij,s}) & \to & \\
  & & \downarrow & & \downarrow & &  \\        
 \stackrel{\alpha_1}{\to} H^{1}(X_{t}) & \stackrel{\alpha_2}{\to} &
 \oplus_{E'_{i,s}} H^{1}(E'_{i,s}) & \stackrel{\alpha_3}{\to} & 
 \oplus_{F'_{ij,s}} H^{1}(F'_{ij,s})& \to &
 H^{2}(X_{t}) \to  \\
  & & \downarrow & & \downarrow & &  \\        
 0 \to \oplus_{F_{ij,s}} H^{0}(F_{ij,s}) & \to &
 \oplus_{E_{i,s} \supset F_{ij,s}}H^{0}(F_{ij, s})  & \to &
 \oplus_{F_{ij,s}} H^{0}(F_{ij,s})  & \to  0\\
 & & \downarrow \beta & & \downarrow & &  \\        
 & & \oplus_{E_{i,s}} H^{2}(E_{i, s}) & \to & 
     \oplus_{F_{ij,s}} H^{2}(F_{ij ,s}) & & 
\end{array}
$$
 Note that the each column is Gysin sequence.
 Then 
 $\dim {\rm Im}(\alpha_1) = 1$ and
 $\dim {\rm Ker}(\alpha_4) = 2n-2$, because
 the singular fibre of $f_{C}$ satisfies the
 properity (1) of Propositon \ref{toroidal_model}.
 Since $\dim {\rm Im}(\alpha_1) + 
 \dim {\rm Im}(\alpha_2) = 2n$,  
 $\dim {\rm Ker}(\alpha_3)= 2n-1$. 
 Hence $\beta$ is not isomorphism. 
 The morphism $\beta$ sends
$$
  H^{0}(F_{ii-1, s}) \oplus H^{0}(F_{ii+1,s}) \ni
   (\lambda_1 , \lambda_2 ) \mapsto 
   (\lambda_1 [F_{ii-1,s}] + \lambda_2 [F_{ii+1,s}])
  \in H^{2}(E_{i,s}).
$$
 If
 $F_{i i-1 ,s}$ and $F_{i i+1 ,s}$ is not
 numerically equivalent, $\beta$ is injection.
 That derives a contradiction and we are done.
\qed

\noindent
 We go back to the proof of Lemma. From Claim \ref{projectivity},
 the nef cone of each irreducible component
 of a singular fibre is generated
 by a section and divisors on an Abelian variety.
 Therefore every big divisor on 
 each irreducible component of singular fibre 
 is ample and projectivity of $f$ can be shown similarly
 as in the case $T = I_{2m}$. \qed

\noindent
  We complete the proof of Proposition \ref{toroidal_model}. \qed

\noindent
 For the using later, we prove the following Lemma.
\begin{lem}\label{flopping}
  Let $f :X \to \Delta^{n}$ be a projective Abelian fibration which
  satisfies the properties (1), (2)
  and (3) of Proposition \ref{toroidal_model}.
  Then every bimeromorphic map $\Phi : X \dasharrow X$
  which commutes with $f$
  can be extended to a biholomorphic morphism, that is,
  $X$ is the unique relative minimal model.
\end{lem} 
\noindent
{\sc Proof .\quad}
       We derive a contradiction assuming that $\Phi$ is not
       isomorphism.
       Since $X$ is a relative minimal model and $\Phi$
       commutes with $f$,
       $\Phi$ is an isomorphism over codimension one
       point of $X$. Let $H$ be an ample divisor on $X$ and
       $H'$ the proper transform by $\Phi$.
       If $\Phi$ is not isomorphism, then $H'$ is not $f$-nef.
       Hence by \cite[Theorem 3-1-1]{KMM},
       there exists an extremal contraction $\nu : X \to \bar{X}$.
       Let $\ell$ be a rational curve which is contracted by $\nu$
       and $E$
       an irreducible component of $f^{-1}(D)$ which
       contains $\ell$.
       Since $f$ satisfies the properties (1) and (2) of Proposition
       \ref{toroidal_model} and
       the discriminat locus $D$
       of $f$ is simply connected, 
       $E \to D$ factors $E \to A \to D$, where
       $E \to A$ is a ${\mathbb P}^{1}$-bundle and
       $A \to D$ is a smooth Abelian fibration.
       Hence $\ell$ must be a fibre of the ruling
       of $E$. Thus $\nu$ contracts every fibre of the ruling
       of $E$ 
       and  
       the exceptional locus of $\nu$ has codimension
       one. 
       That derives a contradiction.        \qed

\begin{defn}\label{definition_of_toroidal_model}
 Let $f : X \to \Delta^{n}$ be a
 projective  Abelian fibration which satisfies
 the properties of Proposition \ref{toroidal_model} (1) - (4).
 We call $X$ as
 a toroidal model of type I or II according to
 the singular fibre of $f$ is an Abelian variety
 or a cycle of several copies of a ${\mathbb P}^{1}$-bundle
 over an Abelian variety.    
\end{defn} 
\noindent
Next we note variations of Hodge structures. 
Let $f :X \to S$ be a projective morphism between smooth
manifold. Assume that there exists an open
set $U$ of $S$ such that $f$ is smooth over $U$ and
$D := S\setminus U$ is a simple normal crossing divisor.
Let $f^{\circ} : f^{-1}(U) \to U$ be the smooth part of $f$.
We denote the upper canonical extension of
${\mathcal H} := R^{i}f^{\circ}_{*}{\mathbb C}$ 
by ${}^{u}{\mathcal H}$ and the lower canonical extension
by ${}^{\ell}{\mathcal H}$.
For the proof of Theorem \ref{main}, we need the following Lemma.
\begin{lem}\label{relation_VHS}
 Let $f :X \to \Delta^{n}$ be a projective morphism from
 a smooth manifold $X$ to a polydisk $\Delta^{n}$.
 Assume that $f$ has the following properties:
\begin{enumerate}
 \item The discriminant locus $D$ of $f$ is smooth and
       $E := (f^{*}D)_{{\rm red}}$ is a simple normal crossing divisor.
 \item The morphism $E \to D$ is a locally trivial deformation.
\end{enumerate}
 Let $f^{*}D = \sum a_i E_i$.
 Then
$$
 {}^{u}{\mathcal F}^{1}({\mathcal H}) 
 \cong f_{*}(\Omega^{1}_{X/\Delta^{n}}(\log E)\otimes \sum (a_i - 1)E_i ).
$$
\end{lem} 
\noindent
{\sc Proof. \quad} If $\dim S = 1$, then this assertion
 is \cite[VIII]{katz}. Thus 
$
 {}^{u}{\mathcal F}^{1}({\mathcal H}) 
$
 and
$
 f_{*}(\Omega^{1}_{X/\Delta^{n}}(\log E)\otimes \sum (a_i - 1)D_i )
$
 are isomorphic in codimension one points of $S$.
 By the assumptions of $f$, \cite[Theorem 3.5]{steenbrink} and
 the argument in the proof of \cite[Lemma 2.11]{torsionfree},
$
 f_{*}(\Omega^{1}_{X/\Delta^{n}}(\log E)\otimes \sum (a_i - 1)D_i )
$ 
 is locally free. 
 Therefore we obtain the assertion of Lemma \ref{relation_VHS}.
\qed

\section{Correspondence between higher direct images and
         cotangent bundles}
 In this section, we prove
 Theorem \ref{main},  Corollary \ref{main_result}
 and Corollary \ref{4dim}. 
 We begin with to consider the relation
 among $R^{i}f_{*}{\mathcal O}_{X}$.
 From the following lemma, it is enough to consider 
 $R^{1}f_{*}{\mathcal O}_{X}$ for the proof of Theorem \ref{main}.
\begin{lem}\label{reduction}
 Let $f :X \to S$ be a Lagrangian fibration
 between smooth projective manifolds. Assume that
 $R^{1}f_{*}{\mathcal O}_{X} \cong \Omega_{S}^{1}$.
 Then 
$$
 R^{i}f_{*}{\mathcal O}_{X} \cong \Omega_{S}^{i}.
$$
\end{lem}
\noindent
{\sc Proof. \quad} Since $\omega_{X} \cong {\mathcal O}_{X}$,
 $R^{i}f_{*}{\mathcal O}_{X}$ are torsion free by 
 \cite[Theorem 2.1]{kollar}.
 Thus $R^{i}f_{*}{\mathcal O}_{X}$ are reflexive sheaves by
 \cite[Corollary 3.9]{torsionfree}.
 Therefore it is enough to show that there exists an open set $U$ of $S$
 such that $R^{i}f_{*}{\mathcal O}_{X}|_{U} \cong \Omega^{i}_{U}$
 and ${\rm codim} (S\setminus U ) \ge 2$.
 Let $D$ be the discriminant locus of $f$.
 We choose an open set $U$ of $S$ such that $D$ is smooth in $U$.
 Since a general fibre of $f$ is an Abelian variety
 ( Liouville's Theorem )
 there exists a morphism
$$\wedge^{i} R^{1}f_{*}{\mathcal O}_{X}|_{U} \stackrel{\alpha_{i}}{\to}
   R^{i}f_{*}{\mathcal O}_{X}|_{U}.
$$
 Then 
 $\alpha$ is injective because
 $R^{1}f_{*}{\mathcal O}_{X}|_{U}$ is torsion free. 
 By \cite[Theorem 2.6]{torsionfree}, 
 $R^{i}f_{*}{\mathcal O}_{X}|_{U} 
 \cong {}^{\ell}{\mathcal Gr}^{i}({\mathcal H})$.
 Considering exponents of the lower extension of
 ${\mathcal H}$,
 $\alpha_{n}$ is not surjective if $\alpha_{i}$ $(0 < i <n)$ is not 
 surjective.
 By \cite[Corollary 7.6]{kollar},
 $R^{n}f_{*}{\mathcal O}_{X} \cong \omega_S$.
 Combining $R^{1}f_{*}{\mathcal O}_{X} \cong \Omega^{1}_{S}$,
 $\alpha_{n}$ is isomorphism. Hence $\alpha_i$
 is isomorphism for $1 \le i \le n$.
\qed

\vspace{3mm}
\noindent
{\sc Remark 4. \quad} The reflexiveness of $R^{i}f_{*}{\mathcal O}_{X}$
  is one of crucial point why we need projectivity.
  This comes from the following decomposition theorem.
\begin{thm}[{\cite[Theorem 3.1]{torsionfree}}]
  Let $f :X \to S$ be a surjective morphism between projective
  variety. If $X$ is smooth, then
$$
  Rf_{*}\omega_{X} \sim_{{\rm q.i.s.}} \sum R^{i}f_{*}\omega_{X}[-i].
$$
\end{thm}
\noindent
 If the above theorem valid for a projective morphism,
 we can replace projectivity of $X$ by
 projectivity of $f$ in the assumptions of Theorem \ref{main}.

\vspace{3mm}

\noindent
 Theorem \ref{main} can be deduced from the following proposition.
\begin{prop}\label{key_proposition}
 Let $f : X \to S$ be a projective Lagrangian fibration with a
 symplectic form $\omega$.
 Then there exists an open set $U$ of $S$ such that
$$
 R^{1}f_{*}{\mathcal O}_{X}|_{U} \cong \Omega_{U}^{1}
$$ 
 and ${\rm codim}(S\setminus U) \ge 2$.
\end{prop}
\noindent
{\sc Proof. \quad}
We prove Proposition \ref{key_proposition} in four steps.

\noindent
{\sc Step 1. \quad}
 First we consider the smooth part of $f$.
 Let $U_0 := S\setminus D$, where $D$
 is the discriminant locus of $f$.
 We consider the following exact sequence:
$$
\begin{array}{ccccccc}
  & & & T_{f^{-1}(U_0)} & \to  & f^{*}T_{U_0}& \to 0 \\
 & & &  \downarrow \omega& & & \\
  0 \to & f^{*}\Omega^{1}_{U_0} & \to & 
      \Omega^{1}_{f^{-1}(U_0)} & \to & 
  \Omega^{1}_{f^{-1}(U_0)/U_0} & \to  0
\end{array}.
$$
 From the above diagram,
 $\omega$ defines a morphism
$$
 f^{*}T_{U_0} \to \Omega^{1}_{f^{-1}(U_0)/U_0 }
$$
 because
 every fibre of $f$ is a Lagrangian subvariety by
 \cite[Theorem 1]{matsu3}. 
 Since $\omega$ is nondegenerate, the above morphism 
 is isomorphism. Therefore we obtain isomorphisms
$$
 T_{U_0} \cong f_{*}\Omega^{1}_{f^{-1}(U_0)/U_0 } 
 \cong {\mathcal F}^{1}({\mathcal H}).
$$
 Taking the dual of above isomorphism,
\begin{equation}\label{smooth_part}
 R^{1}f_{*}{\mathcal O}_{X}|_{U_0} \cong \Omega^{1}_{U_0}.
\end{equation}

\noindent
{\sc Step 2. \quad}
 In order to prove Proposition \ref{key_proposition},
 we extend the isomorphism (\ref{smooth_part})
 over codimension one point of $S$.
 Let $s$ be a point of $D$ such that 
 $S$ and $D$ are smooth around $s$.
 In Step 2, 3 and 4, we will show that the isomorphism (\ref{smooth_part})
 can be extended on a neighborhood of $s$.
 In step 2, we treat the case that $f^{-1}(s)$ satisfies
 the properties of Theorem \ref{proposition_of_classification} (1)
 and $G_1  = G_2 = \{1\}$.
 In this case,
 there exists a polydisk $s \in \Delta^{n}$ of $S$
 and the restriction morphism
 $f : X_{\Delta^{n}} := X\times_{S} \Delta^{n} \to \Delta^{n}$ 
 is a toroidal model of type I or II (See Definition 
 \ref{definition_of_toroidal_model}).
 If $f$ is a toroidal model of type I, then $f$ is smooth
 and we obtain $R^{1}f_{*}{\mathcal O}_{X_{\Delta^{n}}} 
 \cong \Omega^{1}_{\Delta^{n}}$
 around $s$.
 Thus we may assume that $f$ is a toroidal model of type II.
 Let 
$$
 X^{\circ} := \{ x \in X_{\Delta^{n}} | \mbox{$f$ is smooth at $x$}\}.
$$
 We consider the following diagram:
$$
\begin{array}{ccccccc}
   & & & T_{X^{\circ}} & \to  & f^{*}T_{\Delta^{n}}|_{X^{\circ}}& \to 0 \\
 & & &  \downarrow \omega & & & \\
 0 \to & f^{*}\Omega^{1}_{\Delta^{n}}|_{X^{\circ}} & \to & 
      \Omega^{1}_{X^{\circ}} & \to & 
  \Omega^{1}_{X^{\circ}/\Delta^{n}} & \to  0
\end{array}.
$$
 From the above diagram, $\omega$ defines a morphism
 $f^{*}T_{\Delta^{n}}|_{X^{\circ}} \to \Omega^{1}_{X^{\circ}/\Delta^{n}}$
 since every fibre of $f$ is a Lagrangian subvariety by
 \cite[Theorem 1]{matsu3}.
 We investigate the relation between
 $\Omega^{1}_{X^{\circ}/\Delta^{n}}$ and  
 $\Omega^{1}_{X_{\Delta^{n}}/\Delta^{n}}(\log E)|_{X^{\circ}}$. 
 Let us consider the 
 following diagram:
$$
\begin{array}{ccccccc}
  0 \to & f^{*}\Omega^{1}_{\Delta^{n}}|_{X^{\circ}} & \to & 
      \Omega^{1}_{X^{\circ}} & \to & 
  \Omega^{1}_{X^{\circ}/\Delta^{n}} & \to  0\\
  & \downarrow \alpha  & & \downarrow \beta & & \downarrow \gamma & \\
  0 \to & f^{*}\Omega^{1}_{\Delta^{n}}(\log D)|_{X^{\circ}}
   & \to & \Omega^{1}_{X_{\Delta^{n}}}(\log E)|_{X^{\circ}}
   & \to  & \Omega^{1}_{X_{\Delta^{n}}/\Delta^{n}}(\log E) |_{X^{\circ}} 
   & \to 0 .
\end{array}
$$
 We obtain 
 $\Omega^{1}_{X^{\circ}/S}
  \cong \Omega^{1}_{X_{\Delta^{n}}/\Delta^{n}}(\log E)|_{X^{\circ}}$ if
 ${\rm Coker}\alpha \cong {\rm Coker}\beta$. 
 Let $x$ be a point of $f^{-1}(D)$. 
 From Definition \ref{definition_of_toroidal_model}, 
 if we choose a suitable 
 local parameter $z_i$ at $x$ and $t_i$ at $f(x)$, 
 $f$ is written as 
$$
 z_0 z_1 = t_1 , \quad z_i = t_i \quad ( 2 \le i ),
$$
 on a suitable neighborhood of $x$.
 The generater of ${\rm Coker}\alpha$ is $f^{*}dt_1 / t_1$
 and the generaters of ${\rm Coker}\beta$ are $dz_0 / z_0$ and 
 $dz_1 / z_1$.
 The morphsim ${\rm Coker}\alpha \to {\rm Coker}\beta$
 sends 
$$f^{*}dt_1 / t_1 \mapsto dz_0 / z_0 + dz_1 / dz_1 .$$
 If $x \in X^{\circ}$, $z_0 \ne 0$ or $z_1 \ne 0$.
 Hence
 ${\rm Coker}\alpha \to {\rm Coker}\beta$
 is an isomorphsim and we obtain
 $\Omega^{1}_{X^{\circ}/\Delta^{n}} 
 \cong \Omega^{1}_{X_{\Delta^{n}}/\Delta^{n}}
 (\log E)|_{X^{\circ}}$.
 We conclude that
$$
 f^{*}T_{\Delta^{n}} \cong \Omega^{1}_{X_{\Delta^{n}}/\Delta^{n}}(\log E),
$$
 because ${\rm codim}(X\setminus X^{\circ}) \ge 2$.
 Since 
$f_{*}\Omega^{1}_{X_{\Delta^{n}}/\Delta^{n}}
 (\log E) \cong {\mathcal F}^{1}({\mathcal H})$,
 we obtain 
$$R^{1}f_{*}{\mathcal O}_{X_{\Delta^{n}}} 
 \cong \Omega^{1}_{\Delta^{n}}.
$$

\noindent
{\sc Step 3. \quad}
 In this step, we treat that $f^{-1}(s)$
 satisfies the properties of Theorem \ref{proposition_of_classification}
 (1) and $G_1 \ne \{1\}$ or $G_2 \ne \{1\}$.
 By Theorem \ref{proposition_of_classification} (1),
 there exists a following diagram:
$$ 
\begin{array}{ccccccccc}
  X_{\Delta^{n}} 
  & \stackrel{\eta}{\to} & \bar{X}/G_1 & \stackrel{\nu_1}{\leftarrow} 
  & \bar{X} & = & \tilde{X}/G_2 & \stackrel{\nu_2}{\leftarrow} 
  & \tilde{X} \\
  \downarrow & & \downarrow & & \downarrow 
  & & \downarrow & & \downarrow \\
  \Delta^{n} & =  &  \bar{\Delta}^{n}/G_1 & 
  \leftarrow & \bar{\Delta}^{n} & = & \tilde{\Delta}^{n}/G_2 
  & \leftarrow & \tilde{\Delta}^{n}.
\end{array} 
$$
 First we consider the case $G_1 = \{1\}$ and $G_2 \ne \{1\}$.
 Since $\eta$ is isomorphism on
 the smooth locus of $\tilde{X}/G$,
 $X_{\Delta^{n}} \cong \tilde{X}/G_2$. 
 For a symplectic form $\omega$ on $X$,
 $\nu_{2}^{*}\omega$ is
 a symplectic form on $\tilde{X}$
 because $\nu_{2}$ is \'{e}tale. 
 By Step 2,
  $\nu_{2}^{*}\omega$ defines
 an $G_2$-equivariant ismorphism 
 $R^{1}\tilde{f}_{*}{\mathcal O}_{\tilde{X}} \cong 
  \Omega^{1}_{\tilde{\Delta}^{n}}$.
 Take the $G_2$-invariant part, we obtain
$$
 R^{1}f_{*}{\mathcal O}_{X_{\Delta^{n}}} \cong \Omega^{1}_{\Delta^{n}}.
$$
 Next we consider the case $G_1 \ne \{1\}$. 
 Let $U_1$ be the smooth locus of $\bar{X}/G$. Since $\eta$ is
 isomorphism on $U_1$
 and
 ${\rm codim}((\bar{X}/G)\setminus U_1) \ge 2$, 
 we define $\eta_{*}\omega$ as
 the extension of $\omega |_{U_1}$. On the contrary, $\nu_1$ is
 \'{e}tale in codimension one, $\nu_{1}^{*}
 (\eta_{*}\omega)$ is a symplectic
 form on $\bar{X}$. According to the argument
 in the case $G_1 = \{1\}$,
 $\nu_{1}^{*}(\eta_{*}\omega)$ defines a 
 $G_1$-equivariant
 isomorphsim 
$$
 R^{1}f_{*}{\mathcal O}_{\bar{X}} \cong \Omega^{1}_{\bar{\Delta}^{n}}.
$$ 
 Take a $G_1$-invariant part of the above morphism, we obtain 
$$
 R^{1}f_{*}{\mathcal O}_{X_{\Delta^{n}}} \cong \Omega^{1}_{\Delta^{n}}.
$$ 

\noindent
{\sc Step 4. \quad}
 In this step, we treat that $f^{-1}(s)$ 
 satisfies the properties of Theorem \ref{proposition_of_classification}
 (2).
 If
 we choose a suitable neighborfood $\Delta^{n}$
 of $s$, there is an Abelian fibration 
 $\tilde{f} : \tilde{X} \to \tilde{\Delta}^{n}$ and
 the restriction $f : 
 X_{\Delta^{n}} := X\times_{S}{\Delta^{n}} \to \Delta^{n}$ which
 satisfy the following diagram:
$$ 
\begin{array}{ccccc}
  X_{\Delta^{n}} 
& = & \tilde{X}/G & \stackrel{\tau}{\leftarrow}
  & \tilde{X} \\
  \downarrow & & \downarrow & & \downarrow \\
  \Delta^{n} & =  &  \tilde{\Delta}^{n}/G & \leftarrow & 
  \tilde{\Delta}^{n},
\end{array} 
$$
 where $\tau$ is an \'{e}tale morphism.
 Since $\tau$ is \'{e}tale,
 $\tilde{\omega } := \tau^{*}\omega$ is a symplectic form
 where $\omega$ is a symplectic form on $X_{\Delta^{n}}$.
 If $\tilde{\omega}$ defines an isomorphism
 $R^{1}\tilde{f}_{*}{\mathcal O}_{\tilde{X}} \to 
 \Omega^{1}_{\tilde{\Delta}^{n}}$ 
 we obtain a 
 $R^{1}f_{*}{\mathcal O}_{X_{\Delta^{n}}} \cong
 \Omega^{1}_{\Delta^{n}}$ by taking a $G$-invariant part.
 Thus we consider the Abelian fibration $\tilde{f} : \tilde{X} \to 
 \tilde{\Delta}^{n}$.
 By Theorem \ref{proposition_of_classification} (2),
 there is a birational morphism $\nu : W \to \tilde{X}$ such that
 $\tilde{f}\circ \nu : W \to \tilde{\Delta}^{n}$ satisfies
 the assumptions (1) and (2) of Lemma \ref{relation_VHS}.
 Let
 $E := ((\tilde{f} \circ \nu)^{*}\tilde{D})_{red}$ 
 and $(\tilde{f} \circ \nu)^{*}\tilde{D} = \sum a_i E_i$.
 We will show that $\nu^{*}\tilde{\omega}$ defines the morphism
\begin{equation}\label{non_semistable_case}
 (\tilde{f}\circ \nu)^{*}
  T_{\tilde{\Delta}^{n}} \to \Omega^{1}_{W/\tilde{\Delta}^{n}}(\log E)
 \otimes \sum (a_i - 1)E_i .
\end{equation}
 Let $w$ be a point of $(\tilde{f}\circ \nu)^{-1}(\tilde{D})$. 
 On a suitable neighborhood
 of $w$, $\tilde{f}\circ \nu$ is written as
$$
 z^{a_{0}}_{0}\cdots z^{a_{k}}_{k} = t_1 , \quad 
 z_{k+i} = t_i \quad (i \ge 2),
$$
 where $z_i$ are local parameter at $w$ and $t_i$ are
 local parameter at $\tilde{f}\circ\nu (w)$.
 The morphism $T_{W} \to (\tilde{f}\circ \nu)^{*}
 T_{\tilde{\Delta}^{n}}$
 is written as
$$
 \frac{\partial}{\partial z_j } \mapsto 
 a_j z_{0}^{a_{0}}\cdots z_{j}^{a_{j}-1}\cdots z_{k}^{a_k} 
 \frac{\partial}{\partial t_1 }.
$$
 Let $W_0 := W \setminus {\rm Sing}(\tilde{f}\circ \nu)^{-1}(\tilde{D})$.
 If  $w \in W_0$, one of $z_i$ $(0 \le i \le k)$ is nonzero.
 Thus 
$$ (\tilde{f}\circ \nu)^{*}
 T_{\tilde{\Delta}^{n}}\otimes \sum (1 - a_i)E_i |_{W_0}\subset
 {\rm Image}(T_{W}|_{W_0} \to (\tilde{f}\circ \nu)^{*}
 T_{\tilde{\Delta}^{n}}|_{W_0}) .
$$ 
 We cosider the following diagram:  
$$
\begin{array}{ccccccc}
 0   \to & T_{W/\tilde{\Delta}^{n}} & \to
 & T_{W} & \to  & (\tilde{f}\circ \nu)^{*}T_{\tilde{\Delta}^{n}}
 & \\
 & & &  \downarrow  & & & \\
  & (\tilde{f}\circ \nu)^{*}\Omega^{1}_{\tilde{\Delta}^{n}}& \to & 
      \Omega^{1}_{W} & \to & 
  \Omega^{1}_{W/\tilde{\Delta}^{n}} & \to  0 \\
 & \downarrow & &  \downarrow & & \downarrow & \\
  0 \to & (\tilde{f}\circ \nu)^{*}
  \Omega^{1}_{\tilde{\Delta}^{n}}(\log \tilde{D}) & \to & 
  \Omega^{1}_{W}(\log E) & \to & 
  \Omega^{1}_{W/\tilde{\Delta}^{n}}(\log E) & \to 0 .
\end{array}
$$
 By \cite[Theorem 1]{matsu3}, the restriction of $\nu^{*}\tilde{\omega}$
 to every irreducible component of every fibre is identically zero.
 Thus
 $\nu^{*}\tilde{\omega}(\alpha, *) = 0 $ in
 $\Omega^{1}_{W/\tilde{\Delta}^{n}}$
 for
 all elements $\alpha \in T_{W/\tilde{\Delta}^{n}}$. Hence
 $\nu^{*}\tilde{\omega}$ defines the morphism
$$
 (\tilde{f}\circ \nu)^{*}T_{\tilde{\Delta}^{n}}\otimes 
 \sum (1 - a_i )E_i |_{W_0}
   \to 
 \Omega^{1}_{W/\tilde{\Delta}^{n}}(\log E) |_{W_0}.
$$
 Since ${\rm codim}(W\setminus W_0) \ge 2$,
 we obtain the morphsim (\ref{non_semistable_case}).
 Let
$ 
\tilde{X}^{\circ} := \{ \tilde{x} \in \tilde{X} ;
 \mbox{$\tilde{f}$ is smooth at $\tilde{x}$}\}
$
 and $V$ the open set of $\tilde{X}$ such that
 $\nu^{-1}(V) \to V$ is an isomorphism.
 Note that $a_i = 1$ if 
 $E_i \cap \nu^{-1}(X^{\circ}) \ne \emptyset$.
 By similar argument in Step 2, we obtain
 $$ \Omega^{1}_{W/\tilde{\Delta}^{n}}(\log E)\otimes
   (\sum (a_i - 1)E_i) 
   |_{\nu^{-1}(\tilde{X}^{\circ})} \cong
   \Omega^{1}_{\nu^{-1}(\tilde{X}^{\circ})/\tilde{\Delta}^{n}}.
 $$
 Since $\nu^{*}\tilde{\omega}$ is 
 nondegenerate on $\nu^{-1}(V)$,
 the morphsim (\ref{non_semistable_case}) is an isomorphism
 on $\nu^{-1}(\tilde{X}^{\circ} \cap V)$. 
 Thus we obtain an isomorphism
\begin{equation}\label{non_semistable_case_2}
 \tilde{f}^{*}T_{\tilde{\Delta}^{n}}|_{\tilde{X}^{\circ} \cap V} 
 \to \nu_{*} (\Omega^{1}_{W/\tilde{\Delta}^{n}}(\log E)
 \otimes \sum (a_i - 1 )E_i)
 |_{\tilde{X}^{\circ} \cap V}.
\end{equation}
 By Theorem \ref{proposition_of_classification} (2),
 every fibre of $\tilde{f}$ is reduced. Hence
 $X^{\circ} = X\setminus {\rm Sing}(\tilde{f}^{-1}(\tilde{D}))$
 and
 ${\rm codim}(\tilde{X} \setminus \tilde{X}^{\circ} )\ge 2$.
 Combining with ${\rm codim}(\tilde{X}\setminus V) \ge 2$,
 ${\rm codim}(\tilde{X}\setminus (\tilde{X}^{\circ}\cap V)) \ge 2$.
 Since $\tilde{f}^{*}T_{\tilde{\Delta}^{n}}$ is locally free, 
 we obtain
 the morphism (\ref{non_semistable_case_2}) 
 is an isomorphism on $\tilde{X}$ and
$$
 T_{\tilde{\Delta}^{n}} \cong
 \tilde{f}_{*}(\Omega^{1}_{W/\tilde{\Delta}^{n}}(\log E)
 \otimes \sum (a_i - 1)E_i ) .
$$
 By Lemma \ref{relation_VHS},
$\tilde{f}_{*}(\Omega^{1}_{W/\tilde{\Delta}^{n}}
 (\log E)\otimes \sum (a_i - 1)E_i ) 
 \cong {}^{u}{\mathcal F}^{1}({\mathcal H})$. We conclude that
$T_{\tilde{\Delta}^{n}} \cong {}^{u}{\mathcal F}^{1}({\mathcal H})$ and 
$R^{1}\tilde{f}_{*}{\mathcal O}_{\tilde{X}} 
 \cong \Omega^{1}_{\tilde{\Delta}^{n}}$. 
We complete the proof of Proposition \ref{key_proposition}.
\qed

\vspace{2mm}
\noindent
{\sc Proof of Theorem \ref{main}.}
 Under the assumption of Theorem \ref{main},
 $R^{1}f_{*}{\mathcal O}_{X}$ is a reflexive sheaf
 by \cite[Corollary 3.9]{torsionfree}.
 Hence
 we obtain $R^{1}f_{*}{\mathcal O}_{X} \cong \Omega^{1}_{S}$
 by Proposition \ref{key_proposition}.
 Then by Lemma \ref{reduction}, we obtain the assertion of
 Theorem \ref{main}.
\qed

\vspace{2mm}
\noindent
{\sc Proof of Corollary \ref{main_result}. \quad}
 By \cite[Theorem 2]{matsu} and \cite[Theorem 1]{matsu2},
 $f :X \to S$ is a Lagrangian fibration and $S$ is a Fano manifold.
 From
 \cite[Corollary 3.2]{torsionfree},
$$
 h^{k}(X , \omega_{X}) = \sum_{p + q = k} h^{q}(S , R^{p}f_{*}\omega_{X}).
$$
 Since $\omega_{X} \cong {\mathcal O}_{X}$, the right hand side
 of the above equation is the sum of Hodge numbers of $S$
 by Theorem \ref{main}. By Hodge conjugate 
 $h^{k}(X,{\mathcal O}_{X}) = h^{0}(X , \Omega^{k}_{X})$.
 Moreover
$$
 h^{0}(X,\Omega^{k}_{X}) = 
\left\{
\begin{array}{c}
 0  \quad k \equiv 1 \pmod 2\\
 1  \quad k \equiv 0 \pmod 2 
\end{array}
\right.
$$ 
 since $X$ is an irreducible symplectic manifold.
 From the assumption that $S$ is projective, $h^{p,p}(S) \ge 1$.
 Combining the above results,
 we obtain the assertion of Corollary 
 \ref{main_result}.
\qed

\vspace{2mm}
\noindent
{\sc Proof of Corollary \ref{4dim}. \quad}
 We classify $S$ by using the classification of algebraic surfaces.
 We begin with to prove  that $S$ does not contain curves whose
 selfintersection number are negative.
 Assume the contrary.
 Then there is a birational
 morphism $\pi : S \to \bar{S}$ which contracts curves
 to points. The composition morphism $\pi \circ f$ is a 
 Lagrangian fibration and not equidimensional. This
 contradicts to \cite[Theorem 1]{matsu3}.
 Next we investigate the irregularity and Kodaira dimension
 of $S$. From \cite[Theorem 1.1]{addition}, $\kappa (S) \le 0$.
 By 
 similar argument in the proof of Corollary \ref{main_result}, 
 $2q(S) = q(X)$.
 Since $\kappa (X) = 0$, $q(X) \le 4$ and $q(S) \le 2$.
 If $q(S) =2$, $X$ is an Abelian variety.
 From the following Claim,
 $S$ is an Abelian surface.
\begin{claim}\label{Abelian_fibration_criterion}
 Let $X$ be an Abelian variety and $f : X \to S$ be
 a fibre space. Then $S$ is an Abelian variety.
\end{claim} 
\noindent
{\sc Proof. \quad}
 Let $F$ be a general fibre of $f$. Since
 $K_{F}$ is trivial,  $F$ is an Abelian variety by
 \cite[Theorem 12]{Abelian_Characterization}.
 We consider the following exact sequence:
$$
 0 \to F \to X \to X/F \to 0
$$
 Then the quotient morphism is nothing but $X \to S$.
 Thus $X/F \cong S$ and $S$ is an Abelian variety.
\qed

\noindent
 We go back to the proof of Corollary \ref{4dim}.
 Since $S$ does not contain curves whose selfintersection
 number are negative and $\kappa (S) \le 0$,
 $S$ is a ${\mathbb P}^{1}$-bundle
 over an elliptic curve or a hyperelliptic surface
 if $q(S) = 1$.
 Moreover, $S$ is ${\mathbb P}^{2}$ or 
 ${\mathbb P}^{1}\times {\mathbb P}^{1}$ if $q(S) = 0$.
 From these classification, 
 we obtain assertions of Colollary \ref{4dim}.
 \qed

\section{Classification of singular fibres over codimension one points}
In this section, we classify 
singular fibres over codimension one points
of a projective Lagrangian fibration.
\begin{thm}\label{proposition_of_classification}
 Let $f : X \to S$ be a projective Lagrangian fibration and
 $D$ the discriminat locus of $f$. Assume that $\dim X = 2n$,
 $S$ and $D$ is smooth at $s$.
 Then there exists a polydisk $s \in \Delta^{n}$ 
 of $S$ and
 the restriction morphism 
 $f_{\Delta^{n}} : X_{\Delta^{n}} \to \Delta^{n}$ 
 has one of the 
 following two properties:
\begin{enumerate}
 \item There is a toroidal model $\tilde{X}$ of type I or II 
       (See Definition \ref{definition_of_toroidal_model}),
       an Abelian fibration $\bar{X} \to \bar{\Delta}^{n}$
       and
       an action of
       a cyclic group $G_2$ (resp. $G_1$) on $\tilde{X}$
       (resp. $\bar{X}$) which satisfy
       the following diagram.
$$ 
\begin{array}{ccccccccc}
  X_{\Delta^{n}} 
  & \stackrel{\eta}{\to} & \bar{X}/G_1 & \stackrel{\nu_1}{\leftarrow} 
  & \bar{X} & = & \tilde{X}/G_2 & \stackrel{\nu_2}{\leftarrow} 
  & \tilde{X} \\
  \downarrow & & \downarrow & & \downarrow 
  & & \downarrow & & \downarrow \\
  \Delta^{n} & =  &  \bar{\Delta}^{n}/G_1 & 
  \leftarrow & \bar{\Delta}^{n} & = & \tilde{\Delta}^{n}/G_2 
  & \leftarrow & \tilde{\Delta}^{n},
\end{array} 
$$       
       where $\eta$ is isomorphism on the smooth locus
       of $\bar{X}/G_1$ and $\bar{X}/G_1$ is normal.
       Moreover
       $\nu_{1}$ is \'{e}tale
       in codimension one and $\nu_{2}$ is \'{e}tale.
 \item There is an Abelian fibration 
       $\tilde{f} : \tilde{X} \to \tilde{\Delta}^{n}$ and 
       an action of
       a cyclic group $G$ on $\tilde{X}$ 
       which satisfy the following properties:
\begin{enumerate}
 \item $X$ and $\tilde{X}$ satisfy the following diagram:
$$ 
\begin{array}{ccccc}
  X_{\Delta^{n}} 
& = & \tilde{X}/G & \stackrel{\tau}{\leftarrow}
  & \tilde{X} \\
  \downarrow & & \downarrow & & \downarrow \\
  \Delta^{n} & =  &  \tilde{\Delta}^{n}/G & \leftarrow & 
  \tilde{\Delta}^{n},
\end{array} 
$$
 where $\tau$ is an \'{e}tale morphism.
 \item For the discriminant locus
       $\tilde{D}$ of $\tilde{f}$,
       $\tilde{f}^{*} \tilde{D}$ is a reduced divisor.
 \item There is a birational morphsim $\nu : W \to \tilde{X}$
       such that $\tilde{f}\circ \nu : W \to \tilde{\Delta}^{n}$
       satisfies the assumptions (1) and (2)
       of Lemma \ref{relation_VHS}.
\end{enumerate}
\end{enumerate}
\end{thm}
\noindent
{\sc Remark 5. \quad} 
   In the case that $\dim X = 2$, $f$ is a minimal elliptic fibration
   and
   $f$ satisfies the properties of Theorem 
    \ref{proposition_of_classification}(2) if the singular
   fibre of $f$ is a Kodaira singular fibre of type $II$, $III$
   or $IV$. If the singular fibre of $f$ is the other types,
   $f$ satisfies the properties of 
   Theorem \ref{proposition_of_classification}(1).

\vspace{2mm}

\noindent
 Before the proof of Theorem \ref{proposition_of_classification},
 we investigate the monodormy around $s$. 
\begin{lem}\label{monodormy_condition}
 Let $f :X \to S$ be a projective Lagrangian fibration. 
 Assume that the discriminant locus of $f$ 
 and $S$ are smooth at $s$.
 Then there is a polydisk $s \in \Delta^{n}$ of $S$ and the
 restriction $f_{\Delta_{n}}: X_{\Delta^{n}}:= X\times_{S}\Delta^{n} 
 \to \Delta^{n}$ has the following properties:
\begin{enumerate}
 \item
 There is a 
 toroidal model $f^{+} : X^{+} \to \tilde{\Delta}^{n}$ 
 of type I or II
 and the action of
 a cyclic group $G$ on $X^{+}$ which 
 commutes with $f^{+}$ and
 satisfies the
 following diagram:
$$
\begin{array}{ccc}
  X_{\Delta^{n}} 
& \stackrel{\alpha}{\dasharrow} & X^{+}/G \\
  \downarrow & & \downarrow \\
  \Delta^{n} & =  &  \tilde{\Delta}^{n}/G ,
\end{array} 
$$
 where $\alpha$ is a bimeromorphic map.
 \item Let $\tilde{D}$ be the ramification locus of
       $\tilde{\Delta}^{n} \to \Delta^{n}$,
       Then
$$
 \dim F^{1}{\rm Gr}^{W}_{1}H^{1}(
 X^{+}_{\tilde{s}} , {\mathbb C})^{G} \ge n-1 ,
$$
       where $X^{+}_{\tilde{s}}$ be a fibre of 
       $f^{+}$ over a point $\tilde{s}$ of $\tilde{D}$.
 \end{enumerate}
\end{lem}
\noindent
{\sc Proof. \quad} (1)
 By the assumption that the discriminant locus is smooth at $s$,
 we may assume that $S$ and $D$ is smooth. 
 Let $\omega$  be a symplectic form on $X$ and
 $\nu : Y \to X$ a birational morphism such that
 $E = (f\circ \nu )^{*}D$ is a simple normal crossing divisor.
 We denote each component of $E$ by $E_i$.
\begin{claim}\label{restriction}
 Let $E := \sum E_i$ be a simple normal crossing variety.
 Then
$$
 {\rm Gr}^{W}_{k}H^{k}(E,{\mathbb C}) = 
    \{
     (\alpha_i ) \in \oplus H^{k}(E_i ,{\mathbb C}) ;
     \alpha_i |_{E_i \cap E_j } = \alpha_j |_{E_i \cap E_j}
    \}.
$$
\end{claim}
\noindent
{\sc Proof. \quad} 
 Let 
 $$
 E^{[k]} := 
 \coprod_{i_0 < \cdots < i_k} (E_{i_0} \cap \cdots \cap E_{i_k})^{\sim},
 $$
 where $\sim$ means the normalization.
 For an index set $I = \{i_0 , \cdots i_k  \}$, we define 
 an inclusion $\delta^{I}_j$
 $$\delta^{I}_j : E_{i_0} \cap \cdots \cap E_{i_k}
      \to E_{i_0} \cap \cdots \cap E_{i_{j-1}} \cap E_{i_{j+1}} \cap
      \cdots E_{i_k}\quad ( 0 \le j \le k ).$$
 We consider the following spectral sequence \cite[Chapter 4]{GS}:
$$
 E^{p,q}_{1} = H^{q}(E^{[p]}, {\mathbb C}) \Longrightarrow
 E^{p+q} = H^{p+q}(E,{\mathbb C}),
$$
 where $D : E^{p,q}_{1} \to E^{p+1 ,q}_{1}$ is defined by the
 $$
  \bigoplus_{|I| = p} \sum_{j=0}^{p} (-1)^{j}(\delta^{I}_{j})^{*} .
 $$
 Since
 this spectral sequence degenerates at $E_2$ level 
 (\cite[Chapter 4.8]{GS}),
 we deduce
$$
{\rm Gr}^{W}_{k}(H^{k}(E,{\mathbb C}))
 = {\rm Ker}(\oplus_{i} H^{k}(E_i ,{\mathbb C}) \stackrel{D}{\to}
             \oplus_{i<j} H^{k} ((E_i \cap E_j )^{\sim} ,{\mathbb C})).
$$
 Thus we obtain
 the assertion of Claim \ref{restriction} from the definition of $D$.
\qqed

\begin{claim}\label{one_form}
 Let $U$ be an open set of $D$ such that
 $f\circ \nu |_{(f\circ \nu)^{-1}(U)}$ 
 is a locally trivial deformation and 
 every fibre of $f\circ \nu |_{(f\circ \nu)^{-1}(U)}$
 is a simple normal crossing variety. 
 Then every point $s'$ of $U$, 
 $$\dim F^{1}{\rm Gr}^{W}_{1}H^{1}(Y_{s'} , {\mathbb C}) \ge n-1 ,$$
 where $Y_{s'}$ is the fibre of $f\circ \nu $ at $s'$.
\end{claim} 
\noindent
{\sc Proof of Claim. \quad} 
 First
 we show that there is a morphism
\begin{equation}\label{on_form_generaton}
 (f\circ \nu)^{*}T_{U} \to \Omega^{1}_{E_i /U} 
\end{equation}
 for the component
 $E_i \cap (f\circ \nu)^{-1}(U) \ne \emptyset$.
 By the choice of $U$, $E_i \to U$ is a smooth morphism.
 Thus there is a following exact sequence:
$$
 0 \to T_{E_i /U} \to T_{E_i} \to (f\circ \nu)^{*}T_U \to 0.
$$
 From the following diagram,
$$
\begin{array}{ccccccccc}
 0  &\to & T_{E_i} & \to & T_{Y} & \to & N_{E_i / Y} & \to &  0 \\
    &    &         &     &  \downarrow &\nu^{*}\omega  & & & \\
 0 & \to & I_{E_i}/I_{E_i}^{2} & \to & \Omega^{1}_{Y} & \to 
 & \Omega^{1}_{E_i} & \to & 0,
\end{array}
$$
 we define $T_{E_i} \to \Omega^{1}_{E_i}$ by $\nu^{*}\omega$.
 By \cite[Theorem 1]{matsu3}, 
 the restriction $\nu^{*}\omega$ to each component of
 a fibre of $f\circ \nu$ is identically zero. Hence
 $\nu^{*}\omega (\alpha , *) = 0 $ in $\Omega^{1}_{E_i /U}$
 for every element $\alpha  \in T_{E_i /U}$ and
 there is the morphism (\ref{on_form_generaton}).
 We take a point $s' \in U$ and denote local parameters 
 of $U$ at $s'$ by $t_j$, 
 $(2 \le j \le n)$.
 Let  $E_{i, s'  }$ be
 the fibre of $E_i \to U$ over $s'$.
 From the morphism (\ref{on_form_generaton}),
 $\alpha_{ij} := \nu^{*}\omega (\partial /\partial t_j, *)|_{E_{i,s'}}$
 defines an element of $H^{0}(E_{i,s'}, \Omega^{1}_{E_{i,s'}})$.
 We will show that $\alpha_{ij}$, $(2 \le j  \le n)$ are linearly
 independent in  $H^{0}(E_{i,s'}, \Omega^{1}_{E_{i,s'}})$ for
 non $\nu$-exceptional $E_i$.
 Assume the contrary. There exists
 a linear combination 
 $\gamma := \sum \lambda_j \partial /\partial t_j$ such that
 $\gamma \ne 0$ and
 $\nu^{*}\omega(\gamma, *  ) |_{E_{i,s'}}= 0$ in
 $H^{0}(E_{i,s'}, \Omega^{1}_{E_{i,s'}})$.
 From the definition of $U$,
 $E_{i, s'}$ is not contain
 $\nu$-exceptional locus
 if  $E_i$ is not $\nu$-exceptional.
 Thus
 there is an open set
 $V$ of $E_{i,s}$ 
 such that $\nu^{*}\omega$ is nondegenerate over $V$.
 For every point $y \in V$, 
$$
 T_{E_{i,s' },y} = \{ \alpha \in T_{Y,y} 
 ; \nu^{*}\omega (\alpha , \beta) = 0,
 \quad \forall \beta \in T_{E_{i , s' },y} \},
$$
 because $\nu(E_{i,s'})$ is a Lagrangian subvariety.
 However,
 $\gamma  \not\in T_{E_{i,s '},y}$ 
 and 
 $\nu^{*}\omega(\gamma, *  ) |_{E_{i,s'}}= 0$ in
 $H^{0}(E_{i,s'}, \Omega^{1}_{E_{i,s'}})$.
 That derives a contradiction.
 If $E_{i ,s'} \cap E_{i' ,s'} \ne \emptyset$,
 $\alpha_{i , j} = \alpha_{i' ,j}$ on $E_{i,s'} \cap E_{i',s'}$
 by definiton.
 Thus the assertion of Claim \ref{one_form} follows by
 Claim \ref{restriction}.
\qed

\noindent
 We go back to the proof of Lemma \ref{monodormy_condition} (1).
 We choose a neighborhood 
 $\Delta^{n}(t_1 , \cdots , t_n)$ of $s$
 such that $D$ is defined by $t_1 = 0$ and
 $\Delta^{n}\cap U \ne \emptyset$.
 Let $(f\circ \nu)^{*}D = \sum e_i E_i$ and
 $e := {\rm L.C.M.}\{
 e_i; E_i \cap (f\circ \nu)^{-1}(U) \ne \emptyset \}$. 
 We consider a $e$-fold cyclic cover 
 $$d : \tilde{\Delta}^{n} \ni (u, t_2 , \cdots , t_n ) \mapsto 
  (t_1 = u^{e} , t_2 , \cdots , t_n) \in \Delta^{n}$$ 
 and
 the normalization 
 $\hat{Y}$ of 
 $Y \times_{\Delta^{n}}
 \tilde{\Delta}^{n}$. 
 Since $\hat{Y}$ is projective and
 $g^{-1}(u)$ has a reduced component for
 a general point $u$ of $d^{-1}(D)$,
 the Abelian fibration $g : \hat{Y} \to \Delta^{n}$ satisfies
 the assumptions of Lemma \ref{flopping}
 except the monodromy conditions.
 Let  $\tilde{T}$ be the monodromy matrix of
 $R^{1}g^{\circ}_{*}{\mathbb C}$,
 where $g^{\circ}$ is smooth part of $g$.
 We will show that $\tilde{T}$ is unipotent and
 the rank of $\tilde{T}-I$ is
 at most one.
 Let $s'$ be a point of $\Delta^{n}\cap U$.
 We choose
 an unit disk $\Delta^{1}$ in $\Delta^{n}$
 such that $\Delta^{1} \cap D = \{ s'\}$
 and put $\tilde{\Delta}^{1} = d^{-1}(\Delta^{1})$.
 We concentrate the restriction
 morphism 
 $$ g_{\tilde{\Delta}^{1}} : \hat{Y}_{\tilde{\Delta}^{1}}:= 
 \hat{Y}\times_{\tilde{\Delta}^{n}}
 \tilde{\Delta}^{1} \to \tilde{\Delta}^{1}, $$
 because the monodromy matrix of 
 $R^{1}(g_{\tilde{\Delta}^{1}})^{\circ}_{*}{\mathbb C}$ is
 $\tilde{T}$, where $g_{\tilde{\Delta}^{1}}^{\circ}$ is
 smooth part of $g_{\tilde{\Delta}^{1}}$.
 By \cite[Theorem $11^{*}$]{KKMS},
 if we take a suitable resolution $\tilde{Y}_{\tilde{\Delta}^{1}}$ 
 of $\hat{Y}_{\tilde{\Delta}^{1}}$,
 every fibre of
 $\tilde{g}_{\tilde{\Delta}^{1}}: \tilde{Y}_{\tilde{\Delta}^{1}} \to \tilde{\Delta}^{1}$
 is reduced and simple normal crossing.
\begin{claim}\label{injection}
 Let $Z_{\tilde{\Delta}^{1}}$ 
 be a resolution of $\hat{Y}_{\tilde{\Delta}^{1}}$.
 We denote the fibre of  $f\circ \nu$ at $s'$
 by $Y_{s'}$ and
 the fibre of $Z_{\tilde{\Delta}^{1}} \to \tilde{\Delta}^{1}$ at 
 $\tilde{s'} := d^{-1}(s')$ by
 $Z_{\tilde{\Delta}^{1},\tilde{s'}}$ 
 Then the induced morphism 
 $H^{1}(Y_{ s'}, {\mathbb C})
  \to H^{1}(Z_{\tilde{\Delta}^{1}, \tilde{s'}} , {\mathbb C})$
 is injection.
\end{claim}
\noindent
{\sc Proof of Claim. \quad}
 Let $Y_{\Delta^{1}} := Y\times_{\Delta^{n}}\Delta^{1}$.
 Since
 $Y_{\Delta^{1}}$ and $Z_{\tilde{\Delta}^{1}}$
 are deformation retract to 
 $Y_{ s'}$ and $Z_{\tilde{\Delta}^{1}, \tilde{s'}}$
 respectively,
 it is enough to prove that 
 $H^{1}(Y_{\Delta^{1}}, {\mathbb C}) \to 
 H^{1}(Z_{\tilde{\Delta}^{1}}, {\mathbb C})$
 is injective.
 From the isomorphism
 $\hat{Y}_{\tilde{\Delta}^{1}}/G \cong Y_{\Delta^{1}}$, 
 $H^{1}(Y_{\Delta^{1}} , {\mathbb C}) 
 \cong H^{1}(\hat{Y}_{\tilde{\Delta}^{1}}, {\mathbb C})^{G}$
 and $H^{1}(Y_{\Delta^{1}},{\mathbb C}) \to  
 H^{1}(\hat{Y}_{\tilde{\Delta}^{1}},{\mathbb C})$ is injection.
 Moreover $\hat{Y}_{\tilde{\Delta}^{1}}$ is a homology manifold 
 because $\hat{Y}_{\tilde{\Delta}^{1}}$ has only quotient
 singularities.
 (\cite[Proposition 1.4]{steenbrink}).
 Hence
 $H^{1}(\hat{Y}_{\tilde{\Delta}^{1}},{\mathbb C}) 
 \to H^{1}(Z_{\tilde{\Delta}^{1}}, {\mathbb C})$ is injection
 by \cite[Th\`{e}or\'{e}me 8.2.4]{HodgeIII}. 
 \qed

\noindent
 We go back to the proof of Lemma \ref{monodormy_condition}.
 From Claim \ref{injection} and Claim \ref{one_form},
$$
\dim F^{1}{\rm Gr}^{W}_{1}H^{1}
 (\tilde{Y}_{\tilde{s'}}, {\mathbb C}) \ge n-1 .
$$
 Since $\tilde{g}$ is a semistable degeneration,
 $\tilde{T}$ is unipotent.
 We consider the Clemens-Schmid exact sequence.
$$
  0 \to
  {\rm Gr}^{W}_{1}(H^{1}
  (\tilde{Y}_{\tilde{\Delta}^{1}, \tilde{s'}}, {\mathbb C})) \to
  {\rm Gr}^{W}_{1}
  (\tilde{Y}_{\tilde{\Delta}^{1}, \tilde{u}},{\mathbb C})) \stackrel{\tilde{T}-I}{\to}
  {\rm Gr}^{W}_{1}
  (\tilde{Y}_{\tilde{\Delta}^{1}, \tilde{u}},{\mathbb C})), 
$$
 where $\tilde{Y}_{\tilde{\Delta}^{1}, \tilde{u}}$ is a
 general fibre of $\tilde{g}_{\tilde{\Delta}^{1}}$.
 From the above exact sequence, 
 $$\dim {\rm Gr}^{W}_{1}H^{1}
 (\tilde{Y}_{\tilde{\Delta}^{1}, \tilde{u}}, {\mathbb C}) \ge 2n-2 ,$$
 Moreover,
the limit Hodge structure of 
$H^{1}(\tilde{Y}_{\tilde{\Delta}^{1},\tilde{u}},{\mathbb C})$ is
  as follows:
$$
 0 \subset W_0 = {\rm Im}(\tilde{T}-I)
  \subset W_1 = {\rm Ker}(\tilde{T}-I) \subset
  W_2 = H^{1}(\tilde{Y}_{\Delta_1 , \tilde{u}}).
$$
 Since $\dim W_0 = \dim W_2 / W_1$,
 the rank of $\tilde{T}-I$ is at most one.
 Therefore
 $\hat{Y}$ satisfies the assumptions of
 Proposition \ref{toroidal_model} and
 there exists a toroidal model 
 $f^{+} : X^{+} \to \tilde{\Delta}^{n}$  of type I or II
 such that $X^{+}$ is bimeromorphic to $\hat{Y}$.
 The Galois group $G$ of the covering $d$
 acts on $\hat{Y}$ and therefore
 act on $\tilde{X}$ bimeromorphically.
 This action commutes with $f^{+}$. 
 By Lemma \ref{flopping},
 the action of $G$ is biholomorphic and $X^{+}/G$ is
 bimeromorphic to $X_{\Delta^{n}}$ over $\Delta^{n}$.

\noindent
 (2)
 Let $T$ be the monodoromy matrix of 
 $R^{1}(f\circ \nu)^{\circ}_{*}{\mathbb C}$ 
 around $D$ and $\tilde{D} := d^{-1}(D)$.
 For the fibre of $X^{+}_{\tilde{s}}$ of
 $f^{+}$ over $\tilde{s} \in \tilde{D}$,
 the action $G$ on $H^{1}(X^{+}_{\tilde{s}},{\mathbb C})$ 
 is determined by the monodoromy matrix $T$.
 Thus the function 
 $$ \mu (\tilde{s} ) :
 \tilde{s} \mapsto \dim  F^{1}{\rm Gr}^{W}_{1}H^{1}
 (X^{+}_{\tilde{s}} ,{\mathbb C})^{G}$$
 is constant on $\tilde{D}$.
 For the proof of the assertion of Lemma 
 \ref{monodormy_condition} (2),
 it is enough to show $\mu (\tilde{s'}) \ge n-1$
 for $\tilde{s'} = d^{-1}(s')$.
 We consider the
 restriction morphism 
 $X^{+}_{\tilde{\Delta}^{1}}:=
  X^{+}\times_{\tilde{\Delta}^{n}}\tilde{\Delta}^{1} \to 
  \tilde{\Delta}^{1}$.
 Then  
 there is a $G$-equivariant bimeromorphic map
 $\hat{Y}_{\tilde{{\Delta}}^{1}} \dasharrow X^{+}_{\tilde{\Delta}^{1}}$. 
 We take a $G$-equivariant resolution of indeterminacy 
 $Z_{\tilde{\Delta}^{1}}$ of
 $\hat{Y}_{{\tilde{\Delta}}^{1}} \dasharrow X^{+}_{\tilde{\Delta}^{1}}$.
 By Claim \ref{injection},
$ H^{1}
 (Y_{\Delta^{1}, s'} ,{\mathbb C}) \to
  H^{1}(Z_{\tilde{\Delta}^{1}, \tilde{s'}} ,{\mathbb C})
$ 
 is injection.  By Claim \ref{one_form}.
 $\dim F^{1}{\rm Gr}^{W}_{1}H^{1}
 (Z_{\tilde{\Delta}^{1}, \tilde{s'}},{\mathbb C})^{G} \ge n-1 $.
 Because $Z_{\tilde{\Delta}^{1}} \to X^{+}_{\tilde{\Delta}^{1}}$
 is $G$-equivariant bimeromorphic morphism,
 there is a $G$-equivariant isomorphism
 $H^{1}(Z_{\tilde{\Delta}^{1}},{\mathbb C}) \cong
  H^{1}(X^{+}_{\tilde{\Delta}^{1}},{\mathbb C})$.
 Since $Z_{\tilde{\Delta}^{1}}$ and $X^{+}_{\tilde{\Delta}^{1}}$
 are deformation retract to $Z_{\tilde{\Delta}^{1}, \tilde{s'}}$ 
 and $X^{+}_{\tilde{\Delta}^{1} ,\tilde{s'}}$ respectively,
$$
 F^{1}{\rm Gr}^{W}_{1}H^{1}
 (Z_{\tilde{\Delta}^{1}, \tilde{s'}} ,{\mathbb C})^{G}
 \cong  F^{1}{\rm Gr}^{W}_{1}H^{1}
 (X^{+}_{\tilde{s'}} ,{\mathbb C})^{G}.
$$
 Therefore we obtain 
 $\mu (\tilde{s'}) \ge n-1$.
\qed

\noindent
{\sc Proof of Theorem \ref{proposition_of_classification}. \quad}
 Let $f :X \to S$ be 
 a projective Lagrangian fibration.
 By Lemma \ref{monodormy_condition}, there exists a toroidal model
 $f^{+} : X^{+} \to \tilde{\Delta}^{n}$ and an
 action of a cyclic group
 $G$. The assertions of Theorem \ref{proposition_of_classification} 
 follows the following 
 Proposition.
\begin{prop}\label{classification}
 \begin{enumerate}
  \item If $f^{+}$ is a toroidal model of type I, then
        $f^{-1}(s)$ satisfies the properties of 
        Theorem \ref{proposition_of_classification}
        (1) or (2).
   \item If $f^{+}$ is a toroidal model of type II,
         then $f^{-1}(s)$ is 
         satisfies the properties of 
         Theorem \ref{proposition_of_classification} (1).
 \end{enumerate}
\end{prop}
\noindent
{\sc Proof of Proposition \ref{classification}. \quad}
(1)   We begin with to prove that
      the representation $G \to H^{0}(X^{+}_{\tilde{s}}, {\mathbb C})$
      is not trivial, where $X^{+}_{\tilde{s}}$ is the fibre
      of $f^{+}$ at $\tilde{s}$.
      Assume the contrary. Since
      every fibre of $f^{+}$  is an Abelian variety,
      the action of $G$ on $X^{+}$ is fixed point free.
      Thus $X^{+}/G$ is smooth and $K_{X^{+}/G}$ is 
      numerically trivial.
      Since $X^{+}/G$ has no rational curve, $X^{+}/G$ is
      the unique relative minimal model. On the contrary,
      $X_{\Delta^{n}}$ is a relative minimal model. Therefore
      $X_{\Delta^{n}} \cong X^{+}/G$. 
      However, $K_{X^{+}/G} \not\sim 0$ because
      $K_{X^{+}}$ is not $G$-invariant.
      This derives a contradiction, because $K_{X_{\Delta^{n}}} \sim 0$.
      We need the following Lemma
      to prove Proposition \ref{classification} (1).
\begin{lem}\label{equivariant_fibration}
 Let $X^{+}_{\tilde{D}} := (f^{+})^{-1}(\tilde{D})$, where 
 $\tilde{D}$ is the ramification locus
 of $\tilde{\Delta}^{n} \to \Delta^{n}$.
 Then $X^{+}_{\tilde{D}} \to \tilde{D}$ decompose
 a $G$-equivariant
 smooth elliptic fibration  $X^{+}_{\tilde{D}} \to A'$ 
 and a $G$-equivariant
 smooth Abelian fibration $A' \to \tilde{D}$.
 Moreover they
 satisfy the following diagram:
$$
\begin{array}{ccccc}
 X^{+}_{\tilde{D}} &  \to &  A' & \to & \tilde{D} \\
  \downarrow           &     &   \downarrow  & & \downarrow \\
 X^{+}_{\tilde{D}}/G & \to &
   A                   &  \to  &  D,
\end{array}
$$
 where $D$ is the branch locus of $\tilde{\Delta}^{n} \to \Delta^{n}$,
 $A \to D$ are smooth $n-1$ dimensional
 Abelian fibration and $A' /G \cong A$.
\end{lem}
\noindent
{\sc Proof. \quad}
 Let $F$ be a fibre of $X^{+}_{\tilde{D}}\to \tilde{D}$.
 First
 we show that $F/G$ is smooth and $q(F/G) = n-1$.
 Let $T$ be the representation matrix of 
 $\rho : G \to {\rm Aut}H^{0}(F, \Omega^{1}_{F})$.
 Since $G$ is a finite cyclic group and 
 Lemma \ref{monodormy_condition} (2),
 $T = {\rm diag}(\zeta , 1,\cdots ,1)$ 
 under suitable 
 coordinate of $H^{0}(F,\Omega^{1}_{F})$.
 Because $F$ is an Abelian variety,
 the action of $G$ 
 is written $(z_1 , z_2 , \cdots, z_n ) \mapsto
 (\zeta z_1 ,z_2 , \cdots ,z_n)$
 around a fixed point $p$ of $F$,  where $(z_1 , \cdots ,z_n)$
 is a local coordinate of $p$ and
 $\zeta$ is a $n$-th root of unity.
 Note that $\zeta \ne 1$ because $\rho$ is not trivial
 representation. 
 The branch locus of the quotient map $F \to F/G$
 are smooth divisors of $F$ and $F/G$ is smooth.
 Since $\dim H^{0}(F, \Omega^{1}_{F})^{G} = n-1$,
 $q(F/G) = n-1$. 

 Next
 we consider relative Jacobian of $X^{+}_{\tilde{D}}/G \to D$.
 Since $q(F/G) = n-1$, there is a $n-1$-dimensional smooth Abelian  
 fibration $A \to D$ which factor $X^{+}_{\tilde{D}}/G \to D$.
 We prove that $\alpha : X^{+}_{\tilde{D}}/G \to A$ is surjective
 and every fibre is connected. 
 We consider the following diagram
$$
\begin{array}{ccccc}
 X^{+}_{\tilde{D}}  &  \to &  B_1   &  &  \\
  \downarrow           &     &   \downarrow  & &  \\
 X^{+}_{\tilde{D}}/G & \to &
   B_2    &  \to  &  A,
\end{array}
$$
 where $X^{+}_{\tilde{D}}/G \to B_{2}$ is the stein factorization
 of $\alpha$ 
 and $X^{+}_{\tilde{D}} \to B_{1}$ is the stein factorization
 of  $X^{+}_{\tilde{D}} \to B_{2}$.
 Let $A_p$ (resp. $B_{1,p}$, $B_{2,p}$) be the fibre of
 $A \to D$ (resp. $B_1 \to D$, $B_2 \to D$)
 at $p$. 
 Note that the fibre $X^{+}_{p}$ of 
 $X^{+}_{\tilde{D}} \to \tilde{D}$ at $p$
 is an Abelian variety and
 $X^{+}_{p} \to B_{1,p}$ is surjective and connected fibre.
 By Claim \ref{Abelian_fibration_criterion},
 $B_{1,p}$ is an Abelian variety.
 Since $B_{1,p} \to B_{2,p}$ is a finite morphism,
 $\kappa (B_{2,p}) \le \kappa (B_{1,p}) = 0$ by
 \cite[Corollary 9]{Abelian_Characterization}.
 From  \cite[Theorem 13]{Abelian_Characterization},
 there is an \'{e}tale morhpism $\tilde{B}_{2,p} \to B_{2,p}$
 and
 $\tilde{B}_{2,p}$ is isomorphic to the product of
 an Abelian variety and a variety of general type.
 Thus
 $0 \le \kappa (\tilde{B}_{2,p}) = \kappa (B_{2,p})$.
 Therefore
 $\kappa (\tilde{B}_{2,p}) = \kappa (B_{2,p}) = 0$, $\tilde{B}_{2,p}$
 is an Abelian variety and $K_{B_{2,p}}$ is numerically trivial.
 Since 
 $B_{2,p} \to A_p$ is finite and $A_p$ is an
 Abelian variety, $B_{2,p}$ is
 an Abelian variety.
 By the universal properties of relative Jacobian,
 $B_2 \cong A$. The stein factorization
 $X^{+}_{\tilde{D}} \to A'$ 
 of $X^{+}_{\tilde{D}} \to A$ is the desired morphism.
\qqed

\noindent
 We go back to the proof of Proposition \ref{classification} (1).
 By  Lemma \ref{equivariant_fibration},
 there is a $G$-equivariant
 smooth Abelian fibration
 $\pi' : A' \to \tilde{D}$.
 Since $A'/G \cong A$,
 the action of $G$ on $A'$ is a translation
 on each fibre of $\pi'$.
 Let 
 $g$ be a generater of $G$ and
 $m$ the smallest integer such that the action of $g^{m}$ 
 on each fibre of $\pi'$ is trivial.

 First we consider the case that $H = \{ 1\}$. 
 In this case, the action of $G$ is fix point free.
 Hence 
 $X^{+}/G$ is smooth. Moreover $X^{+}/G$
 is the unique relative minimal model over $\Delta^{n}$ since
 it has no rational curves.
 On the contrary, $X_{\Delta^{n}}$ 
 is a relative minimal model over $\Delta^{n}$,
 $X_{\Delta^{n}} \cong X^{+}/G$. 
 Thus $X$ satisfies the properties of Theorem 
 \ref{proposition_of_classification} (1) if we put
 $\tilde{X}= X^{+}$, $G_1 =\{ 1\} $ and $G_2 = G$.

 Next we consider the case that $H \ne \{1\}$.
 Since the action of $H$ on each fibre of $\pi'$ is trivial
 and $X^{+}_{\tilde{D}} \to A'$ is a smooth elliptic fibration,
 every singular locus of $X^{+}/H$
 forms a multisection of 
 $X^{+}_{\tilde{D}}/H \to A'$ and 
 $X^{+}_{\tilde{D}}/H  \to A'$ is a ${\mathbb P}^{1}$-bundle.
 For every singular point $p$ of $X^{+}/H$,
 $(p, X^{+}/H) \cong \mbox{(surface cyclic quotient singularity)}\times 
 ({\mathbb C}^{2n-2},0)$.
 The list of surface cyclic quotient
 singularities which occur above
 is found in \cite[Table 5 (pp 157)]{BPV}.
 By definition, the action of $G/H$ on each fibre of $\pi'$ is non trivial
 translation. Thus the action of $G/H$ on $X^{+}/H$
 is fix point free and the quotient morphism
 $X^{+}/H \to X^{+}/G$ is an \'{e}tale
 morphism. 
 Every singular locus also forms a multisection
 of $X^{+}_{\tilde{D}}/G \to A$ and 
 $X^{+}_{\tilde{D}}/G  \to A$ is a ${\mathbb P}^{1}$-bundle.
 For every singular point  $p$ of $X^{+}/G$, 
 $(p, X^{+}/G) \cong
  \mbox{(surface cyclic quotient singularity)}\times 
 ({\mathbb C}^{2n-2},0)$.
 We will construct
 a suitable resolution $X^{+} /G$ and a relative
 minimal model over $\Delta^{n}$ according to singularities
 of $X^{+}/G$. 


 If every singular point $p$ of $X^{+}/G$,
  $(p, X^{+}/G) \cong \mbox{(Du Val singularity)}\times
  ({\mathbb C}^{2n-2}, 0)$, we take the minimal resolusion 
 $\eta : X' \to X^{+}/G$ by using the minimal resolution
 of Du Val singualrities. Then $X'$
 is a relative minimal model over $\Delta^{n}$.
 The morphism $\eta^{*}(X^{+}_{\tilde{D}}/G) \to A$ 
 is a locally trivial deformation whose fibre is
 a Kodaira singular fibre of type $I^{*}_{0}$, $II^{*}$, 
 $III^{*}$ or $IV^{*}$. Combining that
 $A \to D$ is a smooth
 Abelian fibration, we conclude that
 every irreducible component of
 $\eta^{*}(X^{+}_{\tilde{D}}/G)$
 is a ${\mathbb P}^{1}$-bundle over a smooth Abelian fibration.
 By similar argument in Lemma \ref{flopping},
 $X'$ is 
 the unique minimal model over $\Delta^{n}$.
 Hence $X_{\Delta^{n}} \cong X'$.
 It is verified that $X$ satisfies the properties
 of Theorem \ref{proposition_of_classification} (1)
 by $\tilde{X} = \bar{X} = X^{+}$, $G_1 = G$ and $G_2 = \{1\}$.

 If some singular point $p$ of $X^{+}/G$,
 $(p, X^{+}/G) \not\cong 
 \mbox{(Du Val singularity)} \times
 ({\mathbb C}^{2n-2}, 0)$,
 we consider the minimal resolution $ w : W \to X^{+}/H$,
 which is obtained from the minimal resolution of
 surface quotient singularities.
 Then $w^{*}(X^{+}_{\tilde{D}}/H ) \to A'$ 
 is a locally trivial deformation whose fibre is
 a tree of ${\mathbb P}^{1}$. (cf. \cite[pp 158]{BPV})
 Since $A' \to \tilde{D}/H$ is a smooth Abelian fibration,
 $w^{*}(X^{+}_{\tilde{D}}/H)$ is a simple normal crossing
 divisor and
 the morphism $w^{*}(X^{+}_{\tilde{D}}/H) \to \tilde{D}/H$ 
 is a locally trivial deformation.
 Therefore
 $W$
 satisfies the assumption (1) and (2) of Lemma \ref{relation_VHS}.
 Every irreducible component of 
 $w^{*}(X^{+}_{\tilde{D}}/H)$ is a ${\mathbb P}^{1}$-bundle
 over a smooth Abelian fibration.
 By similar argument in Lemma \ref{flopping},
 $W$ is
 the unique minimal resolution of $X^{+}/H$. Thus
 the action of $G/H$ on $X^{+}/G$ can be lifted on
 $W$ and $W \to X^{+}/H$ is $G/H$-equivariant morphism.
 Hence
 $w^{*}(X^{+}_{\tilde{D}}/H) \to A'$ is
 $G/H$-equivariant.
 By contracting
 ${\mathbb P}^{1}$-bundles along its ruling, 
 (cf. \cite[pp 158]{BPV})
 we obtain a $G/H$-equivariant
 relative minimal model $\tilde{X}$ of $W$ over
 $\tilde{\Delta}^{n}/H$.  
 Let $\tilde{f} : \tilde{X} \to \Delta^{n}/H$ and
 $\tilde{X}_{\tilde{D}/H} := \tilde{f}^{-1}(\tilde{D}/H)$.
 By construction,
 $\tilde{X}$ is smooth and
 every singular fibre of $\tilde{X}$ is reduced.
 The morphism $\tilde{X}_{\tilde{D}/H} \to A'$ 
 is a $G/H$-equivariant locally trivial deformation
 whose fibre is a Kodaira singular fibre of type $II$, $III$ or $IV$. 
 Thus the normalization of every irreducible component
 of $\tilde{X}_{\tilde{D}/H}$ is a ${\mathbb P}^{1}$-bundle
 over a smooth Abelian fibration.
 By similar argument in Lemma \ref{flopping},
 $\tilde{X}$ is the unique minimal model.
 Since
 the action of $G/H$ on $A'$ is non trivial translation,
 the action of $G/H$ on
 $\tilde{X}$ is fixed point free. 
 Thus $\tilde{X} \to \tilde{X}/(G/H)$ is an \'{e}tale morphism and
 $\tilde{X}/(G/H)$ is a relative minimal model. Moreover,
 $\tilde{X}/(G/H)$ is the unique relative minimal model because
 $\tilde{X}$
 is so. Hence 
 $\tilde{X}/(G/H) \cong X_{\Delta^{n}}$. 
 Therefore
 $\tilde{X}$, $X$ and $W$ 
 satisfies the properties 
 of Theorem \ref{proposition_of_classification} (2).
 We finish the proof of Proposition \ref{classification} (1).

\noindent
 (2) First we prove the following Lemma.
\begin{lem}\label{minimal_model_of_H}
  Let
  $\tilde{D}$ be the
  discriminat locus of $\tilde{f}$ and
  $\tilde{f}^{*}\tilde{D} = \sum E_i$.
  We define the subgroup $H$ of $G$ as follows:
$$
  H := \{ g \in G ; \mbox{$g(E_i) = E_i$ for all $i$.}
  \}
$$
  Then there is the unique minimal relative minimal model 
  $\bar{f} : \bar{X} \to \tilde{\Delta}^{n}/H$ 
  of $X^{+}/H \to \tilde{\Delta}^{n}/H$ 
  which have the properties of Proposition
  \ref{toroidal_model} (1), (2) and (3). Moreover
\begin{enumerate}
 \item There is an action of $G/H$ on $\bar{X}$. If an element
       $g$ of $G/H$ preserves every irreducible component of 
       every singular fibre, then $g = 1$. 
 \item For every point of $s$ of the discriminat locus of $\bar{f}$,
       $$
  \dim F^{1}{\rm Gr}^{1}H^{1}(\bar{X}_s , {\mathbb C}) \ge n-1,
       $$      
       where $\bar{X}_{s}$ is the fibre of $\bar{f}$ at $s$.
 \item There exists an \'{e}tale morphism $\nu : \tilde{X} \to \bar{X}$
       and an cyclic group $G'$ such that 
       $\tilde{X}$ is a toroidal model of type II and
       $\tilde{X}/G' \cong \bar{X}$.
\end{enumerate}
\end{lem}

\noindent
{\sc Proof. \quad}
 We begin with investigation of the action of $H$ on
 each $E_i$.
 Let $X^{+}_{\tilde{s}}$ be the fibre of $f^{+}$ over 
 $\tilde{s} \in \tilde{D}$. We denote each irreducible
 component of $X^{+}_{\tilde{s}}$ by $E_{i,\tilde{s}}$ and
 $F_{ij ,\tilde{s}} := E_{i,\tilde{s}} \cap E_{j,\tilde{s}} $. 
 From the following Claim, the action of $H$ on 
 $H^{0}(E_{i,\tilde{s}} , \Omega^{1}_{E_{i,\tilde{s}}})$ and 
 $H^{0}(F_{ij,\tilde{s}}, \Omega^{1}_{F_{ij,\tilde{s}}})$ are trivial.
\begin{claim}\label{normal_crossing_action}
  Let $E =\sum E_{i,\tilde{s}}$ be a cycle of
  several copies of a ${\mathbb P}^{1}$-bundle
  over an Abelian variety such that
  $F_{ij, \tilde{s}} := E_{i,\tilde{s}} 
  \cap E_{j, \tilde{s}}$ forms a section of the ruling
  of $E_{i,\tilde{s}}$ and $E_{j,\tilde{s}}$.
  Assume that there is an action of a cyclic group $G$
  on $E$ such that
  $\dim F^{1}{\rm Gr}^{W}_{1}H^{1}(E, {\mathbb C})^{G} \ge n-1$.
  Then
\begin{enumerate}
 \item If $F_{ij, \tilde{s}}$ is $G$-stable,
       then every element of 
       $H^{0}(F_{ij ,\tilde{s}}, \Omega^{1}_{F_{ij},\tilde{s}})$
       is $G$-invariant.
 \item If $E_{i,\tilde{s}}$ is $G$-stable, then 
       every element of 
       $H^{0}(E_{i, \tilde{s}}, \Omega^{1}_{E_{i},\tilde{s}})$
       is $G$-invariant.
\end{enumerate}
\end{claim} 
\noindent
{\sc Proof. \quad}
 Assume that some $F_{ij,\tilde{s}}$ is $G$-stable and
 there exists a non $G$-invariant element of
 $H^{0}(F_{ij, \tilde{s}}, \Omega^{1}_{F_{ij, \tilde{s}}})$.
 Since every $E_{i,\tilde{s}}$ is a ${\mathbb P}^{1}$-bundle
 over an Abelian variety and every $F_{ij ,\tilde{s}}$ forms
 a section of the ruling,
 there is a $G$-equivariant isomorphism
 $F^{1}{\rm Gr}^{W}_{1}H^{1}(E,{\mathbb C})
 \cong H^{0}(F_{ij,\tilde{s}}, \Omega^{1}_{F_{ij,\tilde{s}}})$ 
 by Claim \ref{restriction}.
 Thus there exists a non $G$-invariant element in
 $F^{1}{\rm Gr}^{W}_{1}H^{1}(E,{\mathbb C})$.
 On the contrary,
   $\dim F^{1}{\rm Gr}^{W}_{1}H^{1}(E, {\mathbb C}) 
   = n-1$. Hence
 $F^{1}{\rm Gr}^{W}_{1}H^{1}(E, {\mathbb C})^{G} 
   = F^{1}{\rm Gr}^{W}_{1}H^{1}(E, {\mathbb C})$. 
 That is a contradiction. If some $E_{i,\tilde{s}}$
 is $G$-stable, we obtain the assertion of Claim by similar 
 argument. \qed

\noindent
 We go back to the proof of Lemma
 \ref{minimal_model_of_H}. 
 Since $F_{ij ,\tilde{s}}$ is an Abelian
 variety, the action of $H$ is translation.
 If the action of $H$ on $F_{ij, \tilde{s}}$ is not trivial, the action
 of $H$ is fixed point free.
 We put $\bar{X} := X^{+}/H$.
 It is easy to check that
 $X^{+}/H$ satisfies the assertions of Lemma.

 We consider the case that
 the action of $H$ on 
 $F_{ij,\tilde{s}}$ is trivial.
  In this case, every point of 
  $F_{ij} := E_{i} \cap E_{j}$ is fixed by $H$.
  We choose a point $p \in F_{ij}$ and
  a local coordinate $z_i$, $(1 \le i \le 2n)$
  at $p$ such that
  $E_{i}$ (resp. $E_{j}$) is defined by $z_1 = 0$. (resp. $z_2 = 0$.) 
  Since $E_{i}$ and $E_{j}$ are 
  stable under the action of $H$, 
  the action of $H$ can be written as
$$
  (z_1 , z_2 , z_3 , \cdots , z_{2n}) \mapsto
  (az_1 , bz_2 , z_3 , \cdots , z_{2n}),
$$
  where $a,b \in {\mathbb C}^{*}$. Thus for
  every singular point $q$ of $X^{+}/H$,
$$
  (q, X^{+}/H) \cong 
  \mbox{(surface cyclic quotient singularity)}\times
  ({\mathbb C}^{2n-2},0).
$$
  The morphism
  $((f^{+})^{*}\tilde{D})/H \to \tilde{D}/H$
  is a locally trivial deformation. 
  Every irreducible component of every fibre is 
  a ${\mathbb P}^{1}$-bundle
  over an Abelian variety.
  The singular locus
  of $X^{+}/H$ is contained in $\cup_{i<j} (F_{ij}/H)$.
  We take the minimal resolution $\nu : Z_{0}  \to X^{+}/H$ 
  by using the minimal resolution
  of a surface cyclic quotient singularity.
  Then the morphism 
  $\nu^{*}(((f^{+})^{*}\tilde{D})/H) \to \tilde{D}/H$
  is a locally trivial deformation.
  Since exceptional locus of
  the minimal resolution of a surface 
  cyclic quotient singularity is a chain of ${\mathbb P}^{1}$,
  every fibre of  $\nu^{*}(((f^{+})^{*}\tilde{D})/H) \to \tilde{D}/H$
  is a cycle of several copies of a ${\mathbb P}^{1}$-bundle
  over an Abelian variety.
  Hence $Z_{0}$ satisfies the properties (1) and (2) of
  Proposition \ref{toroidal_model}.
  Every irreducible component of
  $\nu^{*}(((f^{+})^{*}\tilde{D})/H)$ is
  a ${\mathbb P}^{1}$-bundle over a smooth
  Abelian fibration because $\pi (\tilde{D}/H) = \{1\}$. 
  By similar argument in
  Lemma \ref{flopping},
  $Z_0$ is the unique minimal resolution of $X^{+}/H$.
  Therefore
  the action of $G/H$ can be lifted on $Z_{0}$.
  Let $z_{0} : Z_{0} \to Z_{1}$ be
  a $G/H$-equivariant
  extremal contraction. Let $\ell$ be a curve which
  is contracted by $z_{0}$. We denote the irreducible
  component of $\nu^{*}(((f^{+})^{*}\tilde{D}/H))$ which
  contains $\ell$ by $E_{0}$. Since $E_{0}$ is a 
  ${\mathbb P}^{1}$-bundle over a smooth Abelian fibration,
  $\ell$ is a fibre of the ruling of $E_0$.  
  Thus $z_{0}$ is the contraction
  of ${\mathbb P}^{1}$-bundles along its ruling.
  Hence $Z_{1}$ is smooth and it satisfies
  the properties (1) and (2) of Proposition \ref{toroidal_model}.
  Iterating this process,
  we obtain a $G/H$-equivariant
  minimal model $\bar{X}$. By construction,
  $\bar{X}$
  has the properties (1), (2) and (3) of
  Proposition \ref{toroidal_model}.
  By Lemma \ref{flopping}, $\bar{X}$ 
  is the unique relative minimal model.
  Let $\bar{g}$ be an element of $G/H$.
  By the definition of $H$, $\bar{g}=1$ if
  every irreducible component of singular fibres 
  of $Z_{0} \to \tilde{\Delta}^{n}/H$ are
  stable under
  the action of $\bar{g}$.
  Since  
  at least one $\bar{g}$-orbit is not contracted in the
  contracting process $Z_{0} \to \bar{X}$, $\bar{X}$ 
  satisfies the assertion (1) of Lemma.

  For the proof of (2) and (3) of Lemma, we consider the  
  restriction of $\bar{X}$.
  Let $C$ be a smooth curve of $\tilde{\Delta}^{n}/H$ which
  is stable under the action of $G/H$ and which
  intersects transversally
  with the discriminat locus of $\bar{f}$ at $s$.
  We consider the restrictions $(X^{+}/H)_{C} := 
  X^{+}/H\times_{\tilde{\Delta}^{n}/H}C$, 
  $(Z_{0})_{C} := Z_{0}\times_{\tilde{\Delta}^{n}/H}C$
  and
  $\bar{X}_{C} := \bar{X}\times_{\tilde{\Delta}^{n}/H}C$.
  Since $X^{+}$ satisfies the property (2) of Proposition
  \ref{toroidal_model}, 
  $(X^{+}/H)_{C}$ has only quotient singularities. Hence
  there is an injection 
  $H^{1}((X^{+}/H)_{C}, {\mathbb C}) \to 
   H^{1}((Z_{0})_{C} , {\mathbb C})$. Thus
  $H^{1}((X^{+}/H)_{s}, {\mathbb C}) \to 
   H^{1}((Z_{0})_{s} , {\mathbb C})$ is an injection 
  because $(X^{+}/H)_{C}$ and $(Z_{0})_{C}$ are
  deformation retract to 
  $(X^{+}/H)_{s}$ and $(Z_{0})_{s}$ respectively.
  There is the $G/H$-equivariant bimeromorphic morphism
  $(Z_{0})_{C} \to (X^{+}/H)_{C}$. Moreover
$$  F^{1}{\rm Gr}^{W}_{1}
  H^{1}((X^{+}/H)_{s},{\mathbb C})^{G/H} \cong
   F^{1}{\rm Gr}^{W}_{1}
  H^{1}(X^{+}_{s},{\mathbb C})^{G}.
$$
  We obtain
$$
   \dim F^{1}{\rm Gr}^{W}_{1}
  H^{1}((Z_{0})_{s},{\mathbb C})^{G/H} \ge  n-1,
$$
  from  Lemma \ref{monodormy_condition} (2).
  Since $(Z_{0})_{C} \to \bar{X}_C $ 
  is $G/H$-equivariant bimeromorphic morphism,
  there is a $G/H$-equivariant isomorphism
$
 H^{1}(\bar{X}_C,{\mathbb C}) \to
 H^{1}((Z_{0})_{C},{\mathbb C})
$.
 Hence
$
 F^{1}{\rm Gr}^{W}_{1}
  H^{1}((Z_{0})_{s},{\mathbb C})^{G/H} \cong
 F^{1}{\rm Gr}^{W}_{1}
  H^{1}(\bar{X}_{s},{\mathbb C})^{G/H}
$ 
 because $\bar{X}_{C}$ is deformation retract to $\bar{X}_s$.
 We obtain the assertion (2) of Lemma.
  
  For the
  proof of the assertion (3) of Lemma,
  we need the following Claim.
\begin{claim}\label{multiplicity}
  Let $\bar{f} : \bar{X} \to \tilde{\Delta}^{n}/H$ and
  $\bar{f}^{*}(\tilde{D}/H) = \sum e_i \bar{E}_i$.
  Then $e_i =m$ for all $i$.
\end{claim} 
\noindent
{\sc Proof. \quad} 
  Assume the contrary. Let $\bar{E}_0$ be the
  component such that $e_0$ attains minimum value 
  of $e_i$ and
  $\bar{E}_i$, $(i =1,2)$ next components. We may
  assume that $e_1 > e_0$.
  We consider the restriction $\bar{X}_{C}$ and
  denote the restriction $\bar{E}_i$ to $\bar{X}_{C}$
  by $\bar{E}_{C,i}$.
  Since $\bar{f}$ satisfies the properties (1), (2) and (3)
  of Proposition \ref{toroidal_model},
  each $\bar{E}_{C,i}$ is a ${\mathbb P}^{1}$-bundle over an
  Abelian variety and
 $K_{\bar{E}_{C, 0 }} \sim (K_{\bar{X}_C} + \bar{E}_{C,
 0})|_{\bar{E}_{C, 0}}
 \sim \bar{E}_{C, 0}|_{\bar{E}_{C, 0}}$.
  Let $\ell$ be a fibre of the ruling of $\bar{E}_{C, 0}$.
  Since each intersection $\bar{E}_{C ,i} \cap \bar{E}_{C, j}$ forms
  a section of the ruling of $\bar{E}_{C, i}$. 
$$ 
  (\sum e_i \bar{E}_{C ,i} )\ell =  e_1 + e_2 + e_0 \bar{E}_{C, 0} 
  .\ell = 0.
$$
 Because $\ell$ is a fibre of the ruling of $\bar{E}_{C, 0}$,
 $K_{\bar{E}_{C , 0}}.\ell = -2$. 
 From the assumption $e_1 > e_0$ and $e_2 \ge e_0$,
  this derives a contradiction. \qed

\noindent
  We choose a coordinate $(t_1 ,\cdots , t_n)$ of $\tilde{\Delta}^{n}/H$
  such that $\tilde{D}/H$ is defined by $t_1 = 0$.
  Let $d :\tilde{\Delta}^{n}(s_1 , t_2 ,\cdots , t_n) \to 
  (\tilde{\Delta}^{n}/H)(t_1 , \cdots ,t_n)$ be
  the $m$-fold
  cyclic covering defined by $s^{m}_{1} = t_1$,
  $G'$ the galois group of  $d$   
  and $\tilde{X} :=
  \bar{X}\times_{\tilde{\Delta}^{n}/H}
   \tilde{\Delta}^{n}$. Then
  $\tilde{X} \to \bar{X}$ is \'{e}tale, 
  $\tilde{X}$ is smooth and
  every fibre is reduced. (cf. \cite[Lemma 2.2]{steenbrink})
  We complete the proof of Lemma. \qed

\noindent
  We go back to the proof of Proposition \ref{classification} (2).
  Let us consider the action of $G/H$ on $\bar{X}$, which is 
  constructed in 
  Lemma \ref{minimal_model_of_H}.
  Since the dual graph $\Gamma$ of every singular fibre of $\bar{X}$
  is the Dynkin diagram of type $\tilde{A}_{n}$,
  the action of $G/H$ on $\Gamma$ is  rotation
  or reflection.
  If the action of $G/H$ is rotation, 
  the action of $G/H$ on $\bar{X}$ is fixed point free
  by the assertion (1) of
  Lemma \ref{minimal_model_of_H}.
  The quotient $\bar{X}/(G/H)$ is the unique minimal model
  because $\bar{X}$ is so by Lemma \ref{minimal_model_of_H}.
  Hence $\bar{X}/(G/H) \cong X_{\Delta^{n}}$. 
  Therefore $X$, $\bar{X}$ and $\tilde{X}$
  satisfy the properties of 
  Theorem \ref{proposition_of_classification} (1).

  We consider the case that the action of $G/H$ on $\Gamma$
  is  reflection. If the action of $G/H$ is fixed point free,
  we prove that $\bar{X}/(G/H) \cong X$ and $X$ satisfies the 
  properties (1) of Theorem \ref{proposition_of_classification}
  by similar argument in the case that
  the action of $G/H$ on $\Gamma$ is  rotation.
  Thus we assume that $G/H$ has fixed points.
  We use same notations as in 
  Lemma \ref{minimal_model_of_H} and
  Claim \ref{multiplicity}.
  Since $\bar{X}$ satisfies the properties (1) and (2) of
  Proposition \ref{toroidal_model},
  the morphism $\bar{f}^{*}(\tilde{D}/H) \to \tilde{D}/H$
  is a locally trivial deformation and
  every fibre of $\bar{f}^{*}(\tilde{D}/H) \to \tilde{D}/H$
  is a cycle of several copies of a ${\mathbb P}^{1}$-bundle
  over an Abelian variety. Since $\pi_{1}(\tilde{D}/H) = \{1\}$,
  each $\bar{E}_{i}$ is a ${\mathbb P}^{1}$-bundle
  over a smooth Abelian fibration.
  From the assumption that the action $G/H$ on $\Gamma$ is reflection,
  there are two $G/H$-stable component in $\bar{f}^{*}(\tilde{D}/H)$.
  If $\bar{E}_i$ is $G/H$-stable,
  every one form on $\bar{E}_{i,s}$
  is $G/H$-invariant by 
  Lemma \ref{minimal_model_of_H} (2) and
  Claim \ref{normal_crossing_action}, 
  where $\bar{E}_{i,s}$ is a fibre of $\bar{E}_{i} \to \tilde{D}/H$.
  Since $\bar{E}_{i,s}$ is a ${\mathbb P}^{1}$-bundle
  over an Abelian variety,
  every fibre of the ruling of $\bar{E}_i$ is $G/H$-stable and
  the fixed locus on $\bar{E}_i$ forms
  multisection of the ruling of $\bar{E}_i$.
  If $\bar{F}_{ij} := \bar{E}_{i}\cap \bar{E}_{j}$
  is $G/H$-stable, every one form of 
  $\bar{F}_{ij , s} := \bar{E}_{i ,s}\cap \bar{E}_{j , s}$
  is $G/H$-invariant by 
  Lemma \ref{minimal_model_of_H} (2) and 
  Claim \ref{normal_crossing_action}.
  We obtain that the action of $G/H$ on
  $\bar{F}_{ij}$ is trivial because 
  $\bar{F}_{ij ,s}$ is an Abelian variety and we assume that
  there exists a fixed point.
  Since $\bar{F}_{ij}$ forms a section of the ruling of
  $\bar{E}_i$ and $\bar{E}_j$,
  every connected component of the fixed locus of $\bar{X}$
  forms a multisection of the ruling of some $\bar{E}_i$.
  For every singular point $p$ of $\bar{X}/(G/H)$,
  $(p , \bar{X}/(G/H)) \cong \mbox{($A_1$-singularity)} \times
   ({\mathbb C}^{2n-2}, 0)$.
  We obtain
  a relative minimal model $\bar{f}^{+} : \bar{X}^{+} \to \Delta^{n}$ 
  of $\bar{X}/(G/H)$ over $\tilde{\Delta}^{n}/G = \Delta^{n}$
  by blowing up along singular locus.
  Let $D$ be the discriminat locus of $\bar{f}^{+}$.
  Since every irreducible component
  of $(\bar{f}^{+})^{*}(D)$ is a ${\mathbb P}^{1}$-bundle
  over a smooth Abelian fibration,
  we obtain that
  $\bar{X}^{+}$ is the unique relative minimal model
  by similar argument in Lemma \ref{flopping}.
  Hence
  $X_{\Delta^{n}} \cong \bar{X}^{+}$. From 
  Lemma \ref{minimal_model_of_H},
  there is an \'{e}tale covering $\tilde{X} \to \bar{X}$ such
  that $\tilde{X}$ is a toroidal model of type II. Thus
  $X$ satisfies the properties of Theorem 
  \ref{proposition_of_classification} (1).
\qed

\noindent
 We complete the proof of Theorem \ref{proposition_of_classification}.
\qed

\end{document}